\documentclass[11pt,a4paper]{amsart}
\usepackage{amsmath,amssymb,amsthm}
\usepackage{enumitem}
\usepackage{ifthen,verbatim}
\usepackage[english]{babel}
\usepackage[T1]{fontenc}
\usepackage{floatflt,graphicx,graphics}
\usepackage{a4wide}
\usepackage{cite}
\usepackage{ifthen}

\nonstopmode
\newtheorem{thr}{Theorem}[section]

\newtheorem{cor}{Corollary}[section]

\theoremstyle{definition}
\newtheorem{ex}{Example}[section]
\newtheorem{rem}{Remark}[section]
\def\Re{\mathop{\rm Re}\nolimits}
\def\Im{\mathop{\rm Im}\nolimits}

\numberwithin{figure}{section}
\numberwithin{table}{section}
\usepackage{graphicx}
\usepackage{xcolor}
\begin{document}


\title{One Parameter Families of Conformal Mappings \\ of Bounded Doubly Connected Polygonal Domains}

\author{A.~Dyutin}
\address{Kazan Federal University, Kremlyovskaya str. 35, Tatarstan, 420008,
Russian Federation}
\email{dyutin.andrei2016@yandex.ru}

\author{S.~Nasyrov}
\address{Kazan Federal University, Kremlyovskaya str. 35, Tatarstan, 420008,
Russian Federation}
\email{semen.nasyrov@yandex.ru}


\begin{abstract}
We suggest an approximate method of finding a conformal mapping of an annulus onto an arbitrary bounded doubly connected polygonal  domain. The method is based on the parametric Loewner--Komatu method. We consider smooth one parameter families $\mathcal{F}(z,t)$  of conformal mappings of concentric annuli onto doubly connected  polygonal  domains $\mathcal{D}(t)$ which are obtained from a fixed doubly connected polygonal domain $\mathcal{D}$ by making a finite number of rectilinear slits of variable lengths; meanwhile we do not require the family of domains $\mathcal{D}(t)$ to be monotonous. The integral representation for the conformal mappings $\mathcal{F}(z,t)$ has unknown (accessory) parameters.  We find a PDE for this family and then deduce a system of  PDEs describing dynamics of accessory parameters and conformal modulus of $\mathcal{D}(t)$ when changing the parameter $t$. We note that the right-hand sides of equations in the system of ODEs involve functions standing for speeds of movement of the end points of the slits. This allows us to control completely dynamics of the slits and to seek their agreed change, if we have more than one slit.
We also give some results of numerical calculations which show the efficiency of the suggested method.
\end{abstract}


\keywords{Doubly connected  domain, polygonal domain, conformal modulus, conformal mapping, Schwarz--Christoffel formula, accessory parameters, one parameter family of functions, parametric method, elliptic functions, elliptic integrals, Loewner--Komatu equation.}

\maketitle

\section{Introduction}

In geometric function theory and applications, conformal mappings of polygonal domains play an important role.  In the simply connected case, by the Riemann mapping theorem, there exists a conformal mapping of a given polygonal domain onto either the unit disk or the upper half-plane. The inverse mapping is given by the Schwarz--Christoffel formula (see, e.g. \cite{Driscoll_Tr}). In the doubly connected case, we can map a polygonal domain onto a concentric  annulus (see, e.g.  \cite[ch.~5, \S~1]{Goluzin_book}), and the inverse mapping is expressed via a generalized Schwarz--Christoffel integral, with the help of elliptic and theta functions (see, e.g. \cite{akhiezer,byrd}). First the integral representations for such mappings were obtained by N.~Akhiezer \cite{akhiezer1} (see also \cite[ch.~8, \S.~48]{akhiezer}),  later similar formulas were deduced by G.~M.~Goluzin~\cite{gol1} and Yu.~Komatu \cite{Komatu, komU} (see also \cite[\S~17.5]{Henrici}).

In the multiply connected case, we can map a polygonal domain onto a circular one, i.e. domain bounded with a finite number of disjoint circles  (see, e.g. \cite[ch.~5,
\S~6]{Goluzin_book}).  In 2001 T.~K.~Delillo, A.~R.~Elcrat, J.~A.~Pfaltzgraff \cite{DEP} found a formula for the conformal mapping of the exterior of a finite number of disks onto the exterior of polygonal lines; in 2006 they together with T.~A.~Driscoll  obtained a system of equations to determine unknown parameters in such representations  \cite{DDEP}. D.~G.~Crowdy \cite{C,C1} also obtained a formula for conformal mappings of a multiply connected circular domain onto a bounded and unbounded multiply connected polygonal domain. It is of interest to find conformal mappings of circular polygons onto canonic domains;  we can note that papers \cite{bauer,C2,C3,DK,DK1} with results obtained in this direction.

An undoubted advantage of the Schwarz--Christoffel formula and its generalizations is its simplicity. Its essential disadvantage is the fact that Schwarz--Christoffel integrals contain  unknown parameters. First of all, they are the preimages of the vertices of polygonal lines bounding the domains.  Secondly, the multiply connected case is complicated by the fact that the conformal parameters of domains (so-called conformal moduli) are also unknown. In the doubly connected case, we have only one such parameter which we need to determine, it is called conformal modulus. Now we will recall its definition:
if a doubly connected domain $\mathcal{D}$ is conformally equivalent to an annulus $\{\tau \in \mathbb{C}: q < |\tau| < 1\}$, then its conformal modulus is defined as
$$
\text{Mod} (\mathcal{D}) = \frac{1}{2\pi} \log q^{-1}\,.
$$

There exist various methods of finding unknown (accessory) parameters in Schwarz--Christoffel integrals  and their generalizations. One of the most interesting of them is the Loewner parametric method which has other important applications in the theory of univalent functions (see, e.g., the Schwarz--Christoffel integrals  \cite{bracci}, and the monograph \cite{ppvtof}). Using the method, L.~de Branges \cite{ldb} in 1984 proved the famous Bieberbach conjecture on the estimates of the moduli of the Taylor coefficients for univalent functions in the unit disk \cite{bib}.
We note the recent increase of interest to probabilistic applications of the Loewner method including the stochastic Loewner equation or the Schramm--Loewner equation (SLE) (see, e.g. \cite{gor,bau,dub,gra}).

The parametric Loewner method is based on considering one parameter families of univalent functions depending on a real parameter, that map a canonic domain onto a family of
 domains which are obtained from a fixed one by making an increasing slit, and investigation of a differential equation  which the family satisfy to; it is called the Loewner equation.

The method of finding accessory parameters in Schwarz--Christoffel integrals, based on the  Loewner parametric method, in the simply connected case was first suggested by P.~P.~Kufarev  \cite{kuf1947}; using the Loewner equation he showed that we can determine these parameters by numeric integration of a system of ODEs. In 1960 Yu.~V.~Chistyakov \cite{Chi} first fulfilled Kufarev's method by performing numerical calculations.  In 1979 T.~R.~Hopkins and D.~E.~Roberts \cite{Hopk} also tested this method. Later V.~Ya.~Gutlyanskii and A.~O.~Zaidan \cite{GZ} suggested a modification of Kufarev's method. A slightly different approach, close to one in \cite{GZ}, was suggested by L.~Yu.~Nizamieva \cite{N3}; her study was based on the using of the Hilbert boundary value problem with piecewise constant coefficients. A method of finding accessory parameters in the generalized Schwarz--Christoffel integral, giving a conformal mapping of the half-plane onto a multi-sheeted polygonal domain with branch-points was suggested  in \cite{NK}. In \cite{PN}, the parametric method was developed for the case of several slits. Note that the case of several slits of an arbitrary form was considered in \cite{bom}.

In the multiply connected case, the use of the parametric method is more complicated than in the simply connected one. Making a slit in a non-simply connected domain change its conformal moduli
and this puts the domain with a slit out of the class of domains conformally equivalent to the initial one. Instead, we consider families of conformal mappings, defined on variable domains. In the doubly connected case, there is the Loewner-Komatu equation which is a generalization of the Loewner equation; its development has been started by
Yu.~Komatu \cite{komU} and G.~M.~Goluzin \cite{golO}. Then M.~R.~Kuvaev and P.~P.~Kufarev \cite{kuv,kuvkuf} extended the method to the general case of multiply connected domains.

This paper is devoted to finding accessory parameters in the generalized Schwarz--Christoffel integral with the help of the parametric method in case when we simultaneously make several rectilinear slits of variable lengths in the initial bounded doubly connected domain $\mathcal{D}$; we denote  the corresponding family of polygonal domains with slits by $\mathcal{D}(t)$; here $t$ is a real parameter.

Now we will describe the structure of the paper, the main results and give some comments.

In \S~\ref{sect_ell} we describe some known results from the theory of elliptic functions and elliptic integrals.

In \S~\ref{sect_int}  we give an integral representation for the conformal mapping of a concentric annulus onto a fixed bounded doubly connected polygonal domain  (Theorem~\ref{int_repr}).

Then in \S~\ref{sect_fam} with the use of this representation we consider  smooth one parameter families $\mathcal{F}(z,t)$ of conformal mappings of annuli onto  doubly connected polygonal domains
$\mathcal{D}(t)$. We note that at the same time we make several slits and do not   require that the family $\mathcal{D}(t)$ to be monotonous (in the sense of  inclusion $\subset$).
This allow us to both increase and decrease the lengths of some slits.
We deduce a PDE for the family of such conformal mappings which is a modification of the Loewner--Komatu equation for the case of several rectilinear slits with prescribed dynamics of their endpoints
(Theorem~\ref{family}). From this equation we obtain a system of ODEs  describing the dynamics of accessory parameters and conformal modulus of doubly connected polygonal domains
(Theorem~\ref{accessory}). Solving the Cauchy problem for the system, we find all unknown parameters in the integral representation of
$\mathcal{F}(z,t)$.  We especially note that we can control the location and velocities of end points of the slits,  since their speeds are contained explicitly in the right-hand sides of the obtained differential equations.

We note that we use the approach which was applied earlier for one parameter families of rational functions \cite{NURF2}, conformal mappings of complex tori
\cite{NUT1}, and for doubly connected domains which are the exterior of two rectilinear slits of variable lengths~\cite{DVNslit}, while deducing the system of ODEs. The results obtained in this paper extend the results of \cite{PN}, where one parametric families of polygonal domains with several slits where considered, to the cases of doubly connected domains.

Using the Caratheodory kernel convergence theorem (see, e.g. \cite{Goluzin_book}), we can cut off domains of a more complicated geometry from a given polygonal domain $\mathcal{D}$.  We can also repeat this process and obtain more complicated domains in a few steps. In this case, the method could be fulfilled by consistent solving of several Cauchy problems for systems of ODEs.
The values of accessory parameters and conformal modulus, obtained on every step, except for the final one, give the initial conditions for the  Cauchy problem solved on the next step.

Finally, in \S~\ref{sect_ex} we show the efficiency of the suggested method by giving the results of numeric calculations.
 We take a rectangle with a slit on the axis of symmetry as an initial domain $\mathcal{D}$. To find the initial values of parameters we use the formulas obtained in \S~\ref{sect_int} (Example~\ref{rect_in}).
Then we not only move the rectilinear slit along the axis of symmetry (Example~\ref{ex5.1}), but also cut off polygonal domains moving the interior component of the boundary of the polygonal domain  (Examples~\ref{ex5.2} and \ref{53}). We note that the accuracy of calculations for the preimages of the vertices is about $10^{-9}$ and for the lengths of the sides of boundary polygonal lines is of order~$10^{-6}$.

\section{Preliminary  information}\label{sect_ell}

\textit{Elliptic integrals and elliptic functions}.
First we will recall some facts from the theory of elliptic integrals and elliptic functions (see, e.g. \cite{akhiezer,byrd,lang}).

A meromorphic function in $\mathbb{C}$ is called elliptic, if it has periods $\omega_1$ and $\omega_2$, linearly independent over $\mathbb{R}$. Let us fix an arbitrary point $z_0 \in \mathbb{C}$ and
consider the parallelogram with vertices at the points $z_0$, $z_0 + \omega_1$, $z_0 + \omega_1 + \omega_2$, $z_0 + \omega_2$; it is called the period parallelogram.
We add to this parallelogram all points of the form  $z_0 + r_1 \omega_1 + r_2 \omega_2$ $(0 \le r_1 < 1, 0 \le r_2 < 1)$. For $z_0 = 0$, the period parallelogram is called basic. The order of a nonconstant elliptic functions is the number of its poles in the period parallelogram taking into account multiplicities. Denote by $\Omega =\{ m_1\omega_1 + m_2\omega_2\}$ the lattice generated by the periods $\omega_1$ and $\omega_2$ $(m_1,m_2\in\mathbb{Z})$. Further we will denote by $\omega$ an arbitrary element of the lattice and assume that $\text{Im} (\omega_2/\omega_1)>0$.

One of the main elliptic functions is the Weierstrass $\wp$--function:
$$
\wp(z) = \wp(z;\,\omega_1, \omega_2) = \frac{1}{z^2} + \sum_{\omega\, \in\, \Omega,\, \omega \, \neq \,0}\left(\frac{1}{(z-\omega)^2} - \frac{1}{\omega^2}\right).
$$
It is a meromorphic function in $\mathbb{C}$ with poles of the second order at every point $\omega \in \Omega$. In half-periods the $\wp$--function takes on pairwise different values:
$$
\wp\left(\frac{\omega_1}{2}\right) = e_1, \quad \wp\left(\frac{\omega_2}{2}\right) = e_2, \quad \wp\left(\frac{\omega_1 + \omega_2}{2}\right) = e_3.
$$

We will also need the Weierstrass $\zeta$--function
$$
\zeta(z) = \zeta(z;\, \omega_1, \omega_2) = \frac{1}{z} + \sum_{\omega\, \in\, \Omega,\, \omega\, \neq\, 0}\left(\frac{1}{z-\omega} + \frac{1}{\omega} + \frac{z}{\omega^2}\right).
$$
It has the following properties: $\zeta'(z) = -\wp(z)$ and $\zeta(z + \omega_k) - \zeta(z) = \eta_k$, where $\eta_k = 2\, \zeta({\omega_k}/{2}),\, k = 1, 2$. There is the Legendre relation between periods
$\omega_1$, $\omega_2$ and constants $\eta_1$, $\eta_2$ (see \cite[ch.~3, \S~12]{akhiezer}):
\begin{equation}\label{zetaodnor}
\eta_1\omega_2 - \eta_2\omega_1 = 2\pi{i}\,.
\end{equation}

Another Weierstrass function is
$$
\sigma(z) = \sigma(z;\, \omega_1, \omega_2) = z \prod_{\omega\, \in\, \Omega,\, \omega\, \neq\, 0}\left(1 - \frac{z}{\omega}\right)\exp\left\{{ \frac{z}{\omega} + \frac{z^2}{2\omega^2}}\right\}.
$$
It has the following properties:
$$
\frac{d \log \sigma(z)}{d z} = \zeta(z)\,,\quad \frac{d^2 \log \sigma(z)}{d z^2} = -\wp(z)\,.
$$
$$
\sigma(z \pm \omega_k) = - \exp \left\{{\pm \, \eta_k\left(z \pm \frac{\omega_k}{2}\right)}\right\} \sigma(z)\,,\quad k = 1, 2\,.
$$
We also recall the functions
\begin{equation*}\label{sigma1}
\sigma_1(z) = \pm\, \frac{\sigma(z \pm {\omega_1}/{2})}{\sigma({\omega_1}/{2})}\exp \left\{ { \mp\, \frac{ \eta_1}{ 2}z}\right\} ,
\end{equation*}
\begin{equation*}\label{sigma2}
\sigma_2(z) = \pm\, \frac{\sigma(z \pm {\omega_2}/{2})}{\sigma({\omega_2}/{2})}\exp \left\{ { \mp\, \frac{ \eta_2}{ 2}z}\right\} ,
\end{equation*}
\begin{equation*}\label{sigma3}
\sigma_3(z) = \pm\, \frac{\sigma(z \pm {(\omega_1 + \omega_2)}/{2})}{\sigma({(\omega_1 + \omega_2)}/{2})}\exp \left\{ { \mp\, \frac{ \eta_1 + \eta_2}{ 2}z}\right\} ,
\end{equation*}
and the Weierstrass invariants
$$
g_2 = 60\sum_{\omega\, \in\, \Omega,\, \omega \neq 0}\frac{1}{\omega^4}\,,\quad g_3 = 140\sum_{\omega\, \in\, \Omega,\, \omega\, \neq\, 0}\frac{1}{\omega^6} \,.
$$
Further we will need explicit expressions for the partial derivative of the function $\log \sigma (z) = \log \sigma (z;\, \omega_1,\omega_2)$ with respect to the periods $\omega_1$ and $\omega_2$. In  \cite{NUT1} the following fact is proved.

\begin{thr}\label{der_ln}
{The partial derivatives of the function $ \log \sigma (z) = \log \sigma (z;\, \omega_1,\omega_2)$ with respect to the periods are equal:
$$
\frac{\partial \log \sigma(z)}{\partial \omega_1} = \frac{1}{2\pi i} \left[\frac{1}{2}\, \omega_2 (\wp(z) - \zeta^2(z)) + \eta_2(z \zeta(z) - 1) + \frac{g_2}{24}\, \omega_2  z^2 \right],
$$
\begin{equation}\label{spw2}
\frac{\partial \log \sigma(z)}{\partial \omega_2} = -\frac{1}{2\pi i} \left[\frac{1}{2} \, \omega_1 (\wp(z) - \zeta^2(z)) + \eta_1(z \zeta(z) - 1) -  \frac{g_2}{24}\, \omega_1  z^2 \right].
\end{equation}
}
\end{thr}

We also recall the Jacobi theta-functions, which can be expressed via the Weierstrass elliptic functions (see, e.g. \cite[ch.~IV, \S~19]{akhiezer}):
\begin{equation}\label{theta1sigma}
\vartheta_1(\tau) = \sqrt{\frac{\omega_1}{\pi}} \,(e_1 - e_2)^{{1}/{4}} (e_1 - e_3)^{{ 1}/{ 4}} (e_3 - e_2)^{{ 1}/{ 4}} \exp\left\{ { -\frac{\eta_1 z^2}{ 2\omega_1}}\right\} \sigma(z)\,,
\end{equation}
\begin{equation*}\label{theta2sigma}
\vartheta_2(\tau) = \sqrt{\frac{\omega_1}{\pi}} \, (e_3 - e_2)^{{ 1}/{ 4}}\exp\left\{ { -\frac{ \eta_1 z^2}{ 2\omega_1}}\right\} \sigma_1(z)\,,
\end{equation*}
\begin{equation*}\label{theta3sigma}
\vartheta_3(\tau) = \sqrt{\frac{\omega_1}{\pi}} \, (e_1 - e_2)^{{1}/{4}} \exp\left\{ {-\frac{\eta_1 z^2}{2\omega_1}}\right\} \sigma_3(z)\,,
\end{equation*}
\begin{equation*}\label{theta4sigma}
\vartheta_4(\tau) = \sqrt{\frac{\omega_1}{\pi}} \, (e_1 - e_3)^{{1}/{4}}\exp\left\{ { -\frac{ \eta_1 z^2}{ 2\omega_1}}\right\} \sigma_2(z)\,;
\end{equation*}
here $\tau = z/\omega_1$ and $\exp\{ i\pi \omega_2/\omega_1\} = q$, $|q| < 1$.

We will also need the incomplete elliptic integral of the first kind in the canonic Legendre form:
\begin{equation}\label{FE}
F(z, k) = \int\limits_0^z \frac{d\xi}{\sqrt{(1 - \xi^2)(1 - k^2 \xi^2)}}\,,
\end{equation}
where the parameter $k \in (0,1)$ is called the module of the integral. Putting  $z = 1$ in \eqref{FE}, we obtain
\begin{equation}\label{F1E1}
K(k) = \int\limits_0^1 \frac{d\xi}{\sqrt{(1 - \xi^2)(1 - k^2 \xi^2)}}\,,
\end{equation}
which is called the complete elliptic integral of the first kind. We will also denote $K'(k) = K(k')$ where the parameter $k' = \sqrt{1 - k^2}$ is called the additional module.

The function $z = \textrm{sn}(w,k)$, inverse to that defined by the integral \eqref{FE}, is called the Jacobi elliptic sine; its argument can be expressed via elliptic integrals of the first kind
\eqref{F1E1}.
It is elliptic function of the second kind with the main periods $4K$ and $2iK^{\prime}$.


\section{Integral representation}\label{sect_int}

In this section we will deduce an integral representation for the function mapping conformally an annulus onto a bounded doubly connected polygonal domain. As we noted above, such integral representations was obtained earlier by N.~I.~Akhiezer, G.~M.~Goluzin and Yu.~Komatu.  
Here we give a slightly different representation which will be used further and give its short justification.

Consider a conformal mapping  $w = \mathcal{G}(\tau)$ of an annulus $\mathcal{A} = \{ \tau: q < |\tau| < 1\}$ onto a doubly connected domain $\mathcal{D}$ in $\mathbb{C}$, the exterior boundary component of which is a closed $n_1$-link polygonal line $\Gamma_1$ with vertices at points $w_{1,i_1}$, $1\le i_1\le n_1$, and the interior one is a closed $n_2$-link polygonal line $\Gamma_2$ with vertices at points $w_{2,i_2}$, $1\le i_2\le n_2$.
Denote by $\alpha_{1,i_1}\,\pi$ and $\alpha_{2,i_2}\,\pi$ the interior angles of $\mathcal{D}$ at the vertices $w_{1,i_1}$ and $w_{2,i_2}$. The values of angles satisfy 
\begin{equation}\label{aijnm}
\sum\limits_{i_1 = 1}^{n_1} \alpha_{1,i_1} = n_1 - 2\,,\quad \sum\limits_{i_2 = 1}^{n_2} \alpha_{2,i_2} = n_2 + 2\,.
\end{equation}

Applying the exponential function $\tau(z) = e^{iz}$, we can consider the mapping $\mathcal{F}(z) = \mathcal{G}(e^{iz})$ from the horizontal strip $\{0 < \Im z < \log q^{-1}\}$ onto
$\mathcal{D}$. It maps conformally the rectangle $[0, 2\pi] \times (0, \log q^{-1})$ with  identified vertical sides onto $\mathcal{D}$. We will seek an integral representation of this mapping using the Weierstrass elliptic functions with periods $\omega_1 = 2\pi$ and $\omega_2 = 2 i \log q^{-1}$.

It is evident that $\mathcal{F}$ is analytic in the strip $\{0< \Im z < \log q^{-1}\}$  and is $2\pi$-periodic. Consider the preimages $z_{1,i_1} \in [0,  2\pi]$
and $z_{2,i_2} = x_{2,i_2} + i \log q^{-1} \in [ i \log q^{-1}, 2\pi + i \log q^{-1}]$ of the vertices of $\mathcal{D}$. Denote by $\mathcal{Z}$ the set of points which are obtained from the set of preimages of vertices by all possible shifts to vectors of the lattice $\Omega$, generated by $\omega_1$ and $\omega_2$.
By the Riemann-Schwarz symmetry principle, we extend  $\mathcal{F}$ to $\mathbb{C}\backslash \mathcal{Z}$ as a multiply-valued analytic function. It is easy to see that the function
$\displaystyle \frac{\mathcal{F}''(z)}{\mathcal{F}'(z)}$ is meromorphic and two-periodic in $\mathbb{C}$ with periods $\omega_1$ and $\omega_2$. Consider this function  in the rectangle $\Pi_{\omega} = [0,\omega_1) \times [0, \omega_2/i)$, which is its period parallelogram. We note that the conformal modulus of $\mathcal{D}$ is $\text{Mod} (\mathcal{D}) = \displaystyle \frac{\omega_2}{2\omega_1 i}$. In $\Pi_{\omega}$  the function $\displaystyle \frac{\mathcal{F}''(z)}{\mathcal{F}'(z)}$  has simple poles at the points $z_{1,i_1}$ and $z_{2,i_2}$ with residues $\alpha_{1,i_1}-1$  and $\alpha_{2,i_2}-1$. Therefore, we can restore it by its singularities with the help of the Weierstrass elliptic function $\zeta(z) = \zeta(z;\, \omega_1,\omega_2)$; it has a unique pole in every period parallelogram at the point equivalent to $0$ modulo $\Omega$. We have
$$
\frac{\mathcal{F}''(z)}{\mathcal{F}'(z)} = c  + \sum\limits_{i_1 = 1}^{n_1} (\alpha_{1,i_1} -1)\, \zeta(z - z_{1,i_1}) + \sum\limits_{i_2 = 1}^{n_2} (\alpha_{2,i_2} -1)\, \zeta(z - z_{2,i_2}) \,,
$$
where $c$ is a complex constant. Integrating, we obtain
$$
\log \mathcal{F}' (z) = \log C_1 + c z + \sum\limits_{i_1 = 1}^{n_1} (\alpha_{1,i_1} -1) \log \sigma(z - z_{1,i_1}) + \sum\limits_{i_2 = 1}^{n_2} (\alpha_{2,i_2} -1) \log \sigma(z - z_{2,i_2}) \,.
$$
Then
$$
\mathcal{F}'(z) = C_1 \exp\{c  z\} {\prod\limits_{i_1 = 1}^{n_1} \sigma^{\alpha_{1,i_1} -1}(z - z_{1,i_1})}\,{\prod\limits_{i_2 = 1}^{n_2} \sigma^{\alpha_{2,i_2} -1} (z - z_{2,i_2})} \, ,
$$
consequently,
\begin{equation}\label{F_z_0}
\mathcal{F}(z) = C_1 \int\limits_{0}^{z} \exp\{c \xi\} \prod\limits_{i_1 = 1}^{n_1} \sigma^{\alpha_{1,i_1} -1}(\xi - z_{1,i_1}) \prod\limits_{i_2 = 1}^{n_2} \sigma^{\alpha_{2,i_2} -1} (\xi -
z_{2,i_2}) \, d\xi + C_2\,,
\end{equation}
where $C_1 \neq 0$ and $C_2$ are some complex constants.

Now we will determine the constant $c$ in \eqref{F_z_0}.
The function $\mathcal{F}'(z)$ is periodic with period $\omega_1$. Taking into account the equality
\begin{multline*}
\mathcal{F}^{\prime}(z + \omega_1) = C_1 \exp\{c(z + \omega_1)\} \prod\limits_{i_1 = 1}^{n_1} \sigma^{\alpha_{1,i_1} -1}(z + \omega_1 - z_{1,i_1})
\prod\limits_{i_2 = 1}^{n_2} \sigma^{\alpha_{2,i_2} -1} (z + \omega_1 - (x_{2,i_2} + \omega_2/2)) \\=
\exp\! \Bigg\{\! c \omega_1 + \eta_1 \Bigg(\!\sum\limits_{i_1=1}^{n_1}(\alpha_{1,i_1} -1)(z - z_{1,i_1} + \omega_1/2) \,
\!+ \!\sum\limits_{i_2=1}^{n_2}(\alpha_{2,i_2} -1)(z - x_{2,i_2} + (\omega_1 - \omega_2)/2)\!\Bigg)\!\Bigg\} \mathcal{F}'(z),
\end{multline*}
we conclude that
\begin{multline*} c\omega_1+ \eta_1 \Bigg(\sum\limits_{i_1=1}^{n_1}(\alpha_{1,i_1} -1)(z - z_{1,i_1} + \omega_1/2)\,
\\
+ \sum\limits_{i_2=1}^{n_1}(\alpha_{2,i_2} -1)(z - x_{2,i_2} + (\omega_1 - \omega_2)/2)\Bigg) \equiv  0\quad(\text{mod}\: 2\pi i)\,.
\end{multline*}
Using \eqref{aijnm}, we have
\begin{equation}\label{c_F_z_1}
c\omega_1 - \eta_1\omega_2 - \eta_1 \Bigg(\sum\limits_{i_1=1}^{n_1}(\alpha_{1,i_1} -1)\, z_{1,i_1}    + \sum\limits_{i_2=1}^{n_2}(\alpha_{2,i_2} -1)\, x_{2,i_2} \Bigg) \equiv  0\quad (\text{mod}\;
2\pi i)\,.
\end{equation}
From Legendre's relation \eqref{zetaodnor} and \eqref{c_F_z_1} it follows that
\begin{equation}\label{c_F_z_itog}
c =\frac{\eta_1}{\omega_1}\Bigg(\sum\limits_{i_1=1}^{n_1} (\alpha_{1,i_1} -1)\,z_{1,i_1} + \sum\limits_{i_2=1}^{n_2} (\alpha_{2,i_2} -1)\, x_{2,i_2}\Bigg) + \eta_2\,.
\end{equation}

\begin{thr}\label{int_repr}
{The function $\mathcal{F}$, mapping an annulus $\mathcal{A} = \{\tau:\, q < |\tau| < 1\}$ onto the doubly connected polygonal domain $\mathcal{D}$, has the form \eqref{F_z_0}, where $z = -i \log  \tau$. In \eqref{F_z_0} the points $z_{1,i_1}$ and $z_{2,i_2}$ are the preimages of the vertices of $\mathcal{D}$, the constant $c$ has the form \eqref{c_F_z_itog},  $C_1 \neq 0$ and $C_2$ are some complex constants.
}
\end{thr}

From \eqref{theta1sigma} we can obtain an integral representation of the needed conformal mapping via Jacobi theta-functions.

\begin{cor}\label{repr_theta} {The conformal mapping $\mathcal{F}$, described in Theorem~\ref{int_repr}, has the form
\begin{multline}\label{F_z_teta}
\mathcal{F}(z) = C'_1 \int\limits_{0}^{z} \exp\left\{\Big(-\eta_1 \frac{\omega_2}{\omega_1} + \eta_2\Big)\xi\right\}
\\
\times \prod\limits_{i_1 = 1}^{n_1} \vartheta_1^{\alpha_{1,i_1} -1}\left(\frac{\xi - z_{1,i_1}}{\omega_1}\right) \, \prod\limits_{i_2 = 1}^{n_2} \vartheta_1^{\alpha_{2,i_2} -1}\left(\frac{\xi -
z_{2,i_2}}{\omega_1}\right)  d\xi + C_2\,.
\end{multline}
}
\end{cor}

We note that since theta-functions are expressed by rapidly convergent series, formula \eqref{F_z_teta} is more convenient for numerical calculations than \eqref{F_z_0}.\medskip

Now we will consider an example of conformal mapping onto a doubly connected polygonal domain and find its accessory parameters. This example will be used below in Section~\ref{sect_ex}.

\begin{ex}\label{rect_in}
As a doubly connected polygonal region $\mathcal{D}$, we take the rectangle $\Pi = (-1,1)\times(-b,b)$, $b>0$, with rectilinear slit along the segment $[a_1,a_2]$ of the real axis, $-1 < a_1 < a_2 < 1$. Now we will describe the function
$\mathcal{F}(z)$, that maps conformally the rectangle $[-\pi, \pi] \times (0, \log q^{-1})$ with identified vertical sides onto $\mathcal{D}$.

Since $\mathcal{D}$ is symmetric with respect to the real axis, we may assume that the preimages of the vertices of $\Pi$ are symmetric with respect to the imaginary axis, i.e.
$$
z_{1,1} = \beta_{1,1}, \quad z_{1,2} = \beta_{1,2},\quad  z_{1,3} = -\beta_{1,2}, \quad z_{1,4} = -\beta_{1,1},\quad 0< \beta_{1,1} < \beta_{1,2} < \pi,
$$
besides, $z_{2,1} = \pi + i \log q^{-1}$, $z_{2,2} = i \log q^{-1}$.
In addition, $\mathcal{F}(z_{1,i_1})=w_{1,i_1}$, $1\le i_1\le 4$, where $w_{1,1}= 1 + ib, w_{1,2} = -1 + ib, w_{1,3} = -1 - ib, w_{1,4} = 1 - ib$, and $\mathcal{F}(z_{2,1})=w_{2,1} =a_1$, $\mathcal{F}(z_{2,2})=w_{2,2}=a_2$.

By Theorem~\ref{int_repr}, the function $\mathcal{F}$ that maps an annulus $\mathcal{A} = \{\tau:\;q < |\tau| < 1\}$ onto $\mathcal{D}$, is defined by the formula
\begin{equation}\label{fsigma01}
\mathcal{F}(z) = C_1 \int\limits_{0}^{z} \exp\{c \xi\} \prod\limits_{i_1 = 1}^{4} \sigma^{-1/2}(\xi - z_{1,i_1}) \prod\limits_{i_2 = 1}^{2} \sigma(\xi - z_{2,i_2}) \, d\xi + C_2\,,
\end{equation}
where $z = -i \log \tau$, $c ={\eta_1}/2  + \eta_2$, $\sigma(z) = \sigma(z;\,2\pi,2 i \log q^{-1})$.

Now we will find the accessory parameters in \eqref{fsigma01}. Using the Riemann--Schwarz  symmetry principle, we can consider the conformal mapping $w = \mathcal{G}(\tau)$ of the domain $\mathcal{A}^{+} = \{\tau:\;q < |\tau| <1,\;\text{Im}\,\tau> 0\}$ onto ${{\Pi}}^{+} = \{w:\;-1 < \Re w < 1,\;0 < \Im w < b\}$ and then extend it up to the conformal mapping of $\mathcal{A}$ onto $\mathcal{D}$.

Let $k=k(b)\in(0,1)$ be defined by
$$
k^{1/2} = \frac{\vartheta_2(0, \exp\{-\pi b\})}{\vartheta_3(0, \exp\{-\pi b\})}
$$
and  $c_1 = K(k)$. Then the Jacobi elliptic function $u(w) = \text{sn}(c_1 w, k)$, conformally maps the rectangle ${\Pi}^+$ onto the upper half-plane $\mathbb{H}_u^{+ }$. In this case, the vertices $w_{1,1}^{+} = 1 + ib$, $w_{1,2}^{+} = -1 + ib$, $w_{1,3}^{+} = -1$, $w_{1,4}^{+} = 1$ of the rectangle ${\Pi}^{+}$ will be mapped to the points
$$
u_{1,1} = \frac{1}{k}\,,\quad u_{1,2} = -\frac{1}{k}\,,\quad u_{1,3} = -1,\quad u_{1,4} = 1,
$$
and the points $a_1$ and $a_2$ will be mapped to the points $s_{1} =s_{1}(k)= \text{sn}(c_1 a_1,k)$ and $s_{2} =s_{2}(k)=\text{sn}(c_1 a_2,k)$.

Now consider the function $\zeta(\tau) =\displaystyle i\log \tau + (\pi/2+i\log q^{-1})$, mapping conformally $\mathcal{A}^{+}$ onto the rectangle ${\Pi}_{\mathcal{A}^+} = \{\zeta:\; -\pi/2 < \text{Re}\,\zeta < \pi/2,\; 0 < \text{Im}\,\zeta < \ln q^{-1} \}$. It maps the preimages of the vertices of $\Pi^+$  under the action of the needed function to the points
$$
\zeta_{1,1} = \frac{\pi}{2} - \beta_{1,1} + i \log q^{-1},\ \zeta_{1,2} =\frac{\pi}{2} - \beta_{1,2} + i \log q^{-1}, \ \zeta_{1,3} = -\frac{\pi}{2} + i \log q^{-1}, \ \zeta_{1,4} = \frac{\pi}{2} + i \log q^{-1},
$$
and the points corresponding to the end points of the slits are mapped to the points $\zeta_{2,1} = -\pi/2$ and $\zeta_{2,2} = \pi/2$.

Then we note that the function $v(\zeta) = \text{sn}(c_2 \zeta, \ell)$ maps conformally $\Pi_{\mathcal{A}^+}$ onto the upper half-plane $\mathbb{H}_v^{+}$; here the parameters $c_2$ and $\ell$ satisfy the equalities
\begin{equation} \label{c2ell}
\frac{\pi K^{\prime}(\ell)}{2 K(\ell)}\,=\log q^{-1}, \quad c_2 = \frac{2K(\ell)}{\pi}\,.
\end{equation}
This function maps  the vertices of the rectangle ${\Pi}_{\mathcal{A}^+}$ to the points $v_{1,3} = -1/\ell, v_{1,4} = 1/\ell,  v_{2,1} = -1$, $v_{2,2} = 1$, and the preimages of the vertices of $\Pi$ are mapped to $v_{1,1} = \text{sn}(c_2 \zeta_{1,1}, \ell), v_{1,2} = \text{sn}(c_2 \zeta_{1,2}, \ell)$.
Here the parameter $\ell$  is unknown yet, since it is defined by the value $\log q^{-1}$, connected with the conformal modulus of $\mathcal{D}$ which is also still unknown. Our next problem is to determine the dependence $\ell=\ell(k)$.

First we will construct the Moebius transformation $v = v(u)$ which maps $v_{1,3}$, $v_{1,4}$, $v_{2,1}$ and  $v_{2,2}$ to $u_{1,3}$, $u_{1,4}$, $s_1=s_1(k)$ and $s_2=s_2(k)$. Equating the unharmonic ratios of the corresponding four points, we obtain
$$
\frac{-1/\ell(k) - 1}{-1/\ell(k) + 1} \cdot \frac{1/\ell(k)+1}{1/\ell(k) - 1} = \frac{-1 -s_2(k)}{-1 - s_1(k)} \cdot \frac{1 - s_1(k)}{1 - s_2(k)}\,,
$$
and from this we find
$$
\ell(k) = \frac{s_1(k) s_2(k) - 1 + \sqrt{(s_1^2(k) - 1)(s_2^2(k) - 1)}}{s_1(k) - s_2(k)}\,.
$$

Applying the Moebius transformation $v = v(u)$, we obtain that $v_{{1,1}}=v(u_{{1,1}})$ has the form
$$
v_{{1,1}} = \frac{\left(k(1+\ell) + (1-\ell)\right)s_1 + \left(k(1-\ell) + (1+\ell)\right)s_2 - 2 k s_1 s_2 - 2}{\left(k(\ell + 1) + (\ell - 1)\right)s_1 + \left(k(\ell - 1) + (\ell + 1)\right)s_2 -
2 \ell k s_1 s_2 - 2\ell}\,.
$$
Using \eqref{c2ell} and the fact that the function, inverse to the Jacobi elliptic sine, is expressed via the incomplete elliptic integral, we have
$$
c_2 z_{1,1} = K(\ell) + iK'(\ell) -\int\limits_{1/\ell}^{v_{{1,1}}} \frac{d\xi}{\sqrt{(1 - \xi^2)(1 - \ell^2 \xi^2)}}\,.
$$
On the other side, using the expressions of $c_2$ and $z_{1,1}$ in terms of elliptic functions, we obtain
$$
c_2 z_{1,1} = \frac{2 K(\ell)}{\pi}\left(\frac{\pi}{2} - \beta_{1,1} + i \frac{\pi}{2} \frac{K^{\prime}(\ell)}{K(\ell)}\right),
$$
and from this we find that
\begin{equation*}\label{beta1}
\beta_{1,1}(\ell) =
 \frac{\pi}{2K(\ell)} \int\limits_{1/\ell}^{v_{{1,1}}(\ell)} \frac{d\xi}{\sqrt{(1 - \xi^2)(1 - \ell^2 \xi^2)}}\,.
\end{equation*}

Similarly we deduce:
\begin{equation*}\label{beta2}
\beta_{1,2}(\ell) =
\pi - \frac{\pi}{2K(\ell)} \int\limits_{v_{{1,2}}}^{-1/\ell} \frac{d\xi}{\sqrt{(1 - \xi^2)(1 - \ell^2 \xi^2)}}\,,
\end{equation*}
where
$$
v_{{1,2}} = \frac{\left(k(1+\ell) - (1-\ell)\right)s_1 + \left(k(1-\ell) - (1+\ell)\right)s_2 - 2 k s_1 s_2 + 2}{\left(k(\ell + 1) - (\ell - 1)\right)s_1 + \left(k(\ell - 1) - (\ell + 1)\right)s_2 -
2 \ell k s_1 s_2 + 2\ell}\,.
$$

From \eqref{c2ell} it follows that the conformal modulus of $\mathcal{D}$ is defined by the formula
$$
\text{Mod}(\mathcal{D}) = \frac{K'(\ell)}{4K(\ell)}\,.
$$

In Table~\ref{tab1} we give values of accessory parameters and conformal modulus for some $a_1$ and $a_2$, if $b=0.5$.

\begin{table}[ht]
\caption{Some values of the accessory parameters and the conformal modulus of the rectangle $(-1,1)\times(-0.5,0.5)$, with a linear cut along the segment~$[a_1,a_2]$.}\label{tab1}
{\scriptsize \begin{center}
\begin{tabular}{|{c}|{c}|{c}|{c}|{c}|{c}|}
\hline
$a_1$ & $a_2$ & $\omega_2$ & $z_{1,1}$ & $z_{1,2}$ & $\text{Mod}  (\mathcal{D})$\\
\hline
$-0.75$ & $0.75$  & $1.5034301914585424\,i$ & $0.4710482171084344$ & $2.670544436481146$ & $0.11963917328211081$\\
\hline
$-0.25$ & $0.25$ & $3.3301140313814135\,i$ & $0.19928075807406562$ & $2.9423118955157292$ & $0.26500205457701553$\\
\hline
$-0.1$ & $0.5$  & $2.979450683132465\,i$ & $0.39772860081051964$ & $3.0299504808966304$ & $0.2370971519595281$\\
\hline
$0.1$ & $0.9$ & $2.1183213413271487\,i$ & $1.2511247023650072$ & $3.0991773846339017$ & $0.16857065626463488$\\
\hline
$-0.5$ & $0.5$ & $2.1382753178673837\,i$ & $0.28701506868055354$ & $2.8545775849092414$ & $0.17015854326499394$\\
\hline
\end{tabular}
\end{center}}
\end{table}

\end{ex}

Example~\ref{rect_in} shows that for simple symmetric domains we can find the accessory parameters with the help of elliptic functions. For more complicated domains, we need to apply approximate methods.


\section{Families of conformal mappings} \label{sect_fam}

\begin{figure}[h]
\center{\includegraphics[width=0.82\linewidth]{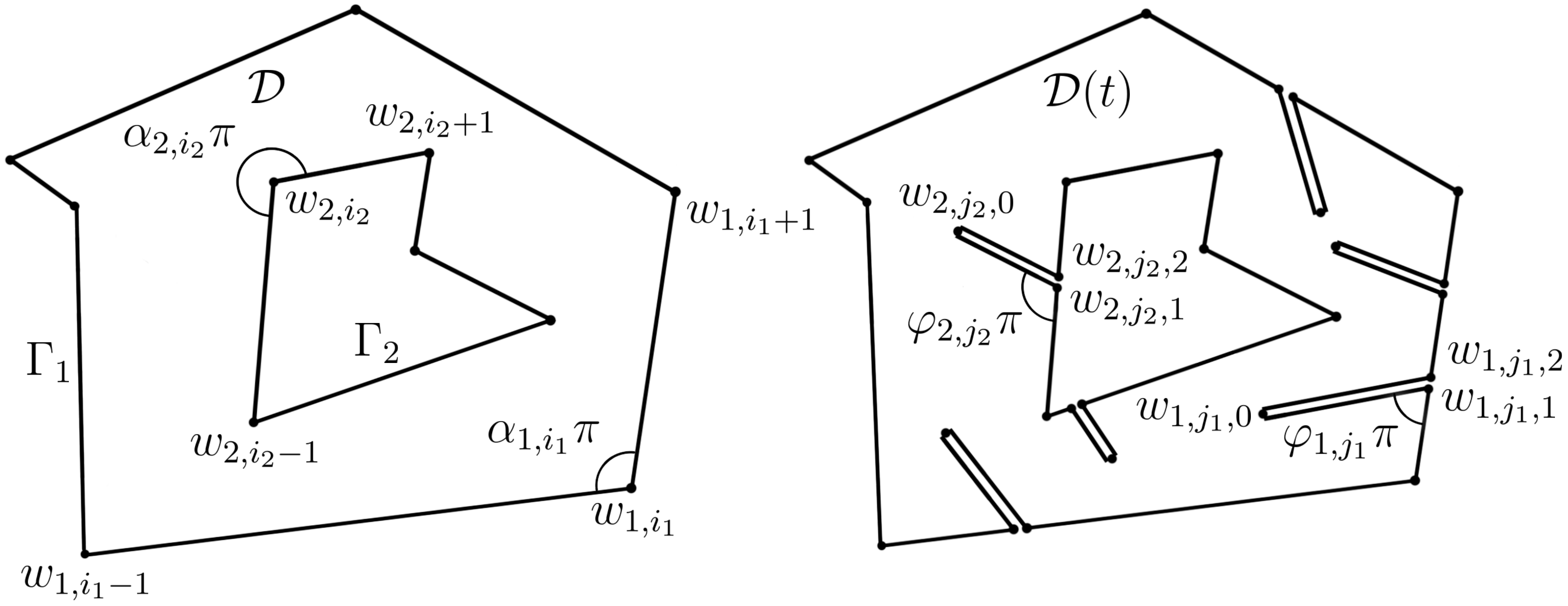}}
\caption{Doubly connected polygonal domain $\mathcal{D}$ and domain $\mathcal{D}(t)$ obtained from $\mathcal{D}$ by drawing a finite number of slits.}\label{ris:image1}
\end{figure}

Now we consider one parameter family of doubly connected domains $\mathcal{D}(t)$, obtained from a fixed doubly connected polygonal domain $\mathcal{D}$ (Fig.~\ref{ris:image1}) by making a few  disjoint rectilinear slits. We will assume that the end points of the slits are smooth functions of the parameter $t$ changing on some segment. For the objects and characteristics, connected with $\mathcal{D}$, we will use the same notation as in Section~\ref{sect_int}.

Let us have $m_1$ slits going from pairwise different points of $\Gamma_1$ at angles $\varphi_{1,j_1} \pi$,  $1\le j_1\le m_1$, and $m_2$ slits from $\Gamma_2$ at angles $\varphi_{2,j_2} \pi$, $1\le j_2\le m_2$ (Fig.~\ref{ris:image1}). For simplicity, we will assume that all the slits proceed from non-angular points of the boundary. The general case can be considered similarly (see Remark~\ref{angle_basepoint} below). Denote by ${z}_{1,j_1,0}(t) \in [0, \omega_1]$ the preimages of the end points of slits proceeding from $\Gamma_1$, and by ${z}_{2,j_2,0}(t) = {x}_{2,j_2,0}(t) + \omega_2(t)/2 \in [\omega_2(t)/2, \omega_1 +\omega_2(t)/2]$ the preimages of the end points of slits proceeding from $\Gamma_2$. We will also denote by $z_{1,j_1,1}(t)$ and $z_{1,j_1,2}(t)$ new angular points of the exterior boundary component with inner angles $\varphi_{1,j_1} \pi$ and $(1 - \varphi_{1,j_1}) \pi$, and by $z_{2,j_2,1}(t)$ and $z_{2,j_2,2}(t)$ we denote new angular points of the interior boundary component with inner angles $\varphi_{2,j_2} \pi$ and $(1 - \varphi_{2,j_2}) \pi$. Now we will describe smooth one parameter families of functions $\mathcal{F}(z,t)$, which conformally map the annuli $\mathcal{A}(t) = \{\tau:\;q(t) < |\tau| < 1\}$ onto $\mathcal{D}(t)$; here $q(t) = \exp \{ - 2 \pi m(t)\}$, $m(t)$ is the conformal modulus of $\mathcal{D}(t)$.

Here,  in contrast to the classical Loewner--Komatu method, we consider the case of several slits and do not assume that the family $\mathcal{D}(t)$ is monotonous with respect to the parameter~$t$.

Now we will find a PDE for the family $\mathcal{F}(z,t)$ and then we will deduce from it a system of ODE to determine the accessory parameters, in particular, the conformal modulus $\text{Mod}(\mathcal{D}(t))$ for all values of parameter~$t$.

So we consider a smooth one parameter family of conformal mappings
\begin{multline}\label{f_sem}
\mathcal{F}(z,t) = C_1(t) \int\limits_{0}^{z} \exp\{c(t)\, \xi \} \prod\limits_{i_1 = 1}^{n_1} \sigma^{\alpha_{1,i_1} -1}(\xi - z_{1,i_1}(t)) \prod\limits_{i_2 = 1}^{n_2} \sigma^{\alpha_{2,i_2} -1}
(\xi - z_{2,i_2}(t))\, \\
\times \prod\limits_{j_1 = 1}^{m_1} s_{1,j_1}(\xi,t) \prod\limits_{j_2 = 1}^{m_2} s_{2,j_2}(\xi,t)\,  d\xi + C_2\,,
\end{multline}
where
$$
s_{1,j_1}(z,t) = \sigma(z - {z}_{1,j_1,0}(t))\, \sigma^{\varphi_{1,j_1} -1}(z - z_{1,j_1,1}(t))\,\sigma^{-\varphi_{1,j_1}} (z - z_{1,j_1,2}(t))\, ,
$$
$$
s_{2,j_2}(z,t) = \sigma(z - {z}_{2,j_2,0}(t))\, \sigma^{\varphi_{2,j_2} - 1}(z - z_{2,j_2,1}(t))\,\sigma^{-\varphi_{2,j_2}}(z - z_{2,j_2,2}(t))\,  ,
$$
\begin{multline*}
c(t) = \frac{\eta_1(t)}{\omega_1}\Bigg[\sum\limits_{i_1=1}^{n_1} (\alpha_{1,i_1} -1)\,z_{1,i_1}(t)+\sum\limits_{i_2=1}^{n_2} (\alpha_{2,i_2} -1)\,x_{2,i_2}(t)\, \\
+ \sum\limits_{j_1=1}^{m_1} \Big({z}_{1,j_1,0}(t) + (\varphi_{1,j_1}-1)\,z_{1,j_1,1}(t) -\varphi_{1,j_1} \, z_{1,j_1,2}(t)\Big) \\
+  \sum\limits_{j_2=1}^{m_2}\Big({x}_{2,j_2,0}(t) + (\varphi_{2,j_2} - 1)\, x_{2,j_2,1}(t) - \varphi_{2,j_2}\,x_{2,j_2,2}(t) \Big) \Bigg] + \eta_2(t)\,.
\end{multline*}
Here $\sigma(z) = \sigma(z;\, \omega_1,\omega_2(t))$, $\, \eta_1(t) = 2 \zeta({\omega_1}/{2};\,\omega_1,\omega_2(t))\,$ and $\, \eta_2(t) = 2\zeta({\omega_2(t)}/{2}; \,\omega_1,$ $\omega_2(t))$.
The conformal modulus of $\mathcal{D}(t)$ equals
\begin{equation*}\label{mod}
\text{Mod}(\mathcal{D}(t)) = \displaystyle  \frac{\omega_2(t)}{2\omega_1 i}\,.
\end{equation*}

For a fixed $t$, the function $\mathcal{F}(z, t)$ is periodic with period $\omega_1 = 2\pi$. By the symmetry principle, we extend $\mathcal{F}(z, t)$ through the segments with endpoints which are the preimages of the vertices of the boundary components. It is evident that the extension satisfies the equalities
\begin{equation*}\label{f_z_t_w1_it}
\mathcal{F}(z + \omega_1,t) = \mathcal{F}(z, t)\,,\quad \mathcal{F}(z + \omega_2(t),t) = e^{ i\alpha }\mathcal{F}(z, t)+\beta \,.
\end{equation*}

Differentiating with respect to $t$ and $z$, we obtain
$$
\dot{\mathcal{F}}(z + \omega_1,t) = \dot{\mathcal{F}}(z,t),\,\quad \dot{\omega}_2(t)\, \mathcal{F}'(z + \omega_2(t),t) + \dot{\mathcal{F}}(z + \omega_2(t),t) = e^{ i\alpha} \dot{\mathcal{F}}(z,
t)\,,
$$
$$
\mathcal{F}'(z + \omega_1,t) = \mathcal{F}'(z, t)\,,\quad \mathcal{F}'(z + \omega_2(t),t) = e^{ i\alpha} \mathcal{F}'(z, t)\,.
$$
Here and further the dot above a letter denotes the  differentiation with respect to $t$, and the prime denotes differentiation with respect to $z$. Therefore, we have
$$
\frac{\dot{\mathcal{F}}\left(z + \omega_1,t\right)}{\mathcal{F}'(z + \omega_1,t)} = \frac{\dot{\mathcal{F}}(z, t)}{\mathcal{F}'(z, t)}\,,
$$
$$
\frac{\dot{\mathcal{F}}\left(z + \omega_2(t),t\right)}{\mathcal{F}'(z + \omega_2(t),t)} + \dot{\omega}_2(t) = \frac{\dot{\mathcal{F}}(z, t)}{\mathcal{F}'(z, t)}\,.
$$
Consequently, the function
$$
\mathcal{H}(z, t) := \frac{\dot{\mathcal{F}}(z, t)}{\mathcal{F}'(z, t)}
$$
satisfies
\begin{equation}\label{h_z_t_w1_w2_it}
\mathcal{H}\left(z + \omega_1, t\right) - \mathcal{H}(z, t) = 0\,,\quad \mathcal{H}\left(z + \omega_2(t), t\right) - \mathcal{H}(z, t) = - \dot{\omega}_2(t)\,.
\end{equation}

Now we will write the Taylor expansion of $\mathcal{F}(z, t)$ near the point $z_{1,k_1}(t)$:
$$
\mathcal{F}(z, t) = A_0 + A_1(t)\left(z - z_{1,k_1}(t)\right)^{\alpha_{1,k_1}} + \ldots \,,\quad 1 \le k_1 \le n_1 \,;
$$
here $A_0$ does not depend on~$t$.
Then
$$
\mathcal{F}'(z, t) = \alpha_{1,k_1}\, A_1(t) \left(z - z_{1,k_1}(t)\right)^{\alpha_{1,k_1} -1} + \ldots \,,
$$
$$
\dot{\mathcal{F}}(z, t) = - \alpha_{1,k_1}\, A_1(t) \left(z - z_{1,k_1}(t)\right)^{\alpha_{1,k_1} -1} \dot{z}_{1,k_1}(t) +  \ldots  \,
$$
and, therefore,
\begin{equation}\label{z_1k_p_t_it}
\dot{z}_{1,k_1}(t) = -\mathcal{H}(z_{1,k_1}(t), t)\,, \quad 1 \le k_1 \le n_1\,.
\end{equation}
Similarly, we obtain
\begin{equation}\label{z_2k_p_t_it}
\dot{z}_{2,k_2}(t) = -\mathcal{H}(z_{2,k_2}(t), t)\,, \quad 1 \le k_2 \le n_2\,,
\end{equation}
\begin{equation}\label{z_1k_p_t_it1}
\dot{z}_{1,\ell_1,j}(t) = -\mathcal{H}(z_{1,\ell_1,j}(t), t)\,, \quad 1 \le \ell_1 \le m_1\,,\quad j =1,2\,,
\end{equation}
\begin{equation}\label{z_2k_p_t_it1}
\dot{z}_{2,\ell_2,j}(t) = -\mathcal{H}(z_{2,\ell_2,j}(t), t)\,, \quad 1 \le \ell_2 \le m_2\,, \quad j =1,2\,.
\end{equation}

Now we will write the Taylor expansion for $\mathcal{F}(z, t)$ near the points ${z}_{1,\ell_1,0}(t)$:
\begin{equation}\label{1l1}
\mathcal{F}(z, t) = \mathcal{E}_{1,\ell_1}(t) + \frac{1}{2}\,{\mathcal{F}''\left({z}_{1,\ell_1,0}(t),t\right)}\left(z - {z}_{1,\ell_1,0}(t)\right)^2 +\ldots\,,\quad 1 \le \ell_1 \le m_1\,,
\end{equation}
where $\mathcal{E}_{1,\ell_1}(t) = \mathcal{F}(z_{1,\ell_1,0}(t),t)$ is the end point of the $\ell_1$th slit, entering from $\Gamma_1$. From the expansion it follows that
\begin{equation*}\label{f_z_t_prime01}
\mathcal{F}'(z, t) = \mathcal{F}''\left({z}_{1,\ell_1,0}(t),t\right) \left(z - {z}_{1,\ell_1,0}(t)\right) +\ldots\,,
\end{equation*}
$$
\dot{\mathcal{F}}(z, t) = \dot{\mathcal{E}}_{1,\ell_1} (t) -  \mathcal{F}''\left({z}_{1,\ell_1,0}(t),t\right) \left(z - {z}_{1,\ell_1,0}(t)\right) \dot{{z}}_{1,\ell_1,0}(t) +\ldots\,,
$$
and, thus,
$$
\mathcal{H}(z, t) = \frac{\mathcal{L}_{1,\ell_1}(t)}{z - {z}_{1,\ell_1,0}(t)} + \text{O}(1)\,,\quad z \rightarrow {z}_{1,\ell_1,0}(t)\,,
$$
where
$$
\mathcal{L}_{1,\ell_1}(t) := \frac{\dot{\mathcal{E}}_{1,\ell_1}(t)}{\mathcal{F}''\left({z}_{1,\ell_1,0}(t),t\right)}, \quad 1 \le \ell_1 \le  m_1\,.
$$

Similarly,
$$
\mathcal{L}_{2,\ell_2}(t) := \frac{\dot{\mathcal{E}}_{2,\ell_2}(t)}{\mathcal{F}''\left({z}_{2,\ell_2,0}(t),t\right)}, \quad 1 \le \ell_2 \le m_2\,.
$$

Now we will find $\mathcal{F}''\left({z}_{1,\ell_1,0}(t),t\right)$. We have
\begin{multline*}
\mathcal{F}''(z, t) = C_1(t) \exp\{c(t)\, z\} \prod\limits_{i_1 = 1}^{n_1} \sigma^{\alpha_{1,i_1} -1}(z - z_{1,i_1}(t)) \prod\limits_{i_2 = 1}^{n_2} \sigma^{\alpha_{2,i_2} -1} (z - z_{2,i_2}(t))\,\\
\times \prod\limits_{j_1 = 1}^{m_1} s_{2,j_1}(z, t) \prod\limits_{j_2 = 1}^{m_2} s_{2,j_2}(z, t) \Bigg[\mathcal{Q}(z,t) + \sum\limits_{j_1 = 1}^{m_1} \mathcal{Q}_{1,i_1}(z,t) + \sum\limits_{j_2 =
1}^{m_2} \mathcal{Q}_{2,i_2}(z,t) \Bigg],
\end{multline*}
where
$$
\mathcal{Q}(z,t) = c(t) + \sum\limits_{i_1 = 1}^{n_1} (\alpha_{1,i_1} - 1)\, \zeta(z - z_{1,i_1}(t)) + \sum\limits_{i_2 = 1}^{n_2} (\alpha_{2,i_2} -1)\, \zeta(z - z_{2,i_2}(t))\, ,
$$
$$
\mathcal{Q}_{1,i_1}(z,t) = \zeta (z - {z}_{1,j_1,0}(t)) + (\varphi_{1,j_1} -1)\, \zeta(z - z_{1,j_1,1}(t)) - \varphi_{1,j_1} \, \zeta (z - z_{1,j_1,2}(t)) \, ,
$$
$$
\mathcal{Q}_{2,i_2}(z,t) = \zeta (z - {z}_{2,j_2,0}(t)) + (\varphi_{2,j_2} -1)\, \zeta(z - z_{2,j_2,1}(t)) - \varphi_{2,j_2} \, \zeta (z - z_{2,j_2,2}(t)) \, .
$$
Consequently, as $z \rightarrow {z}_{1,\ell_1,0}(t)$, we obtain
\begin{multline}\label{f_z_t_pp_z1k2} \mathcal{F}''\left({z}_{1,\ell_1,0}(t),t\right) = C_1(t) \exp\{c(t)\, {z}_{1,\ell_1,0}(t)\}\, \\
\times \prod\limits_{i_1 = 1}^{n_1} \sigma^{\alpha_{1,i_1} -1}({z}_{1,\ell_1,0}(t) - z_{1,i_1}(t))\, \prod\limits_{i_2 = 1}^{n_2} \sigma^{\alpha_{2,i_2} -1} ({z}_{1,\ell_1,0}(t) - z_{2,i_2}(t))\,\\
\times  \prod\limits_{j_1 = 1}^{m_1} \sigma^{\varphi_{1,j_1} -1}({z}_{1,\ell_1,0}(t) - z_{1,j_1,1}(t))\, \sigma^{-\varphi_{1,j_1}} ({z}_{1,\ell_1,0}(t) - z_{1,j_1,2}(t))\, \prod\limits_{j_2 =
1}^{m_2} s_{2,j_2}({z}_{1,\ell_1,0}(t),t)\,.
\end{multline}
Similarly,
\begin{multline}\label{f_z_t_pp_z2k2}
\mathcal{F}''\left({z}_{2,\ell_2,0}(t),t\right) = C_1(t) \exp\{c(t)\, {z}_{2,\ell_2,0}(t)\} \, \\
\times \prod\limits_{i_1 = 1}^{n_1} \sigma^{\alpha_{1,i_1} -1}({z}_{2,\ell_2,0}(t) - z_{1,i_1}(t))\, \prod\limits_{i_2 = 1}^{n_2} \sigma^{\alpha_{2,i_2} -1} ({z}_{2,\ell_2,0}(t) - z_{2,i_2}(t))\,\\
\times \prod\limits_{j_1 = 1}^{m_1} s_{1,j_1}({z}_{2,\ell_2,0}(t),t) \prod\limits_{j_2 = 1}^{m_2} \sigma^{\varphi_{2,j_2} -1}({z}_{2,\ell_2,0}(t) - z_{2,j_2,1}(t))\, \sigma^{-\varphi_{2,j_2}}
({z}_{2,\ell_2,0}(t) - z_{2,j_2,2}(t))\,.
\end{multline}

The function
$$
\mathcal{S}(z, t) := \mathcal{H}(z, t) - \sum\limits_{j_1 = 1}^{m_1} \mathcal{L}_{1,j_1}(t)\, \zeta(z - {z}_{1,j_1,0}(t)) - \sum\limits_{j_2 = 1}^{m_2} \mathcal{L}_{2,j_2}(t)\, \zeta(z -
{z}_{2,j_2,0}(t))
$$
has removable singularities only  at the points which are the preimages of the end points of the slits and at points equivalent to them modulo the lattice. Consequently, it has a holomorphic extension to the whole plane $\mathbb{C}$. From \eqref{h_z_t_w1_w2_it} we have
$$
\mathcal{S}(z + \omega_1, t) - \mathcal{S}(z, t) = -\sum\limits_{j_1 = 1}^{m_1} \mathcal{L}_{1,j_1}(t)\, \zeta(z + \omega_1 - {z}_{1,j_1,0}(t)) -  \sum\limits_{j_2 = 1}^{m_2} \mathcal{L}_{2,j_2}(t)\,
\zeta(z + \omega_1 - {z}_{2,j_2,0}(t)) \,
$$
$$
+ \sum\limits_{j_1 = 1}^{m_1} \mathcal{L}_{1,j_1}(t)\, \zeta(z - {z}_{1,j_1,0}(t)) + \sum\limits_{j_2 = 1}^{m_2} \mathcal{L}_{2,j_2}(t)\, \zeta(z - {z}_{2,j_2,0}(t))
$$
\begin{equation}\label{S_z+w1}
=  - \eta_1(t)\Bigg( \sum\limits_{j_1 = 1}^{m_1} \mathcal{L}_{1,j_1}(t) + \sum\limits_{j_2 = 1}^{m_2} \mathcal{L}_{2,j_2}(t)\Bigg),
\end{equation}
\begin{equation}\label{S_z+w2}
\mathcal{S}(z + \omega_2(t), t) - \mathcal{S}(z, t) = -\dot{\omega}_2(t) - \eta_2(t)\Bigg( \sum\limits_{j_1 = 1}^{m_1} \mathcal{L}_{1,j_1}(t) + \sum\limits_{j_2 = 1}^{m_2}
\mathcal{L}_{2,j_2}(t)\Bigg).
\end{equation}
Because of \eqref{S_z+w1} and \eqref{S_z+w2}, the function $\mathcal{S}$ grows at infinity no faster that a linear one, therefore, $\mathcal{S}(z,t) = \alpha(t)z + \beta(t)$. Thus,
\begin{equation}\label{Hab}
\mathcal{H}(z, t) = \sum\limits_{j_1 = 1}^{m_1} \mathcal{L}_{1,j_1}(t)\, \zeta(z - {z}_{1,j_1,0}(t)) +\sum\limits_{j_2 = 1}^{m_2} \mathcal{L}_{2,j_2}(t)\, \zeta(z - {z}_{2,j_2,0}(t))  +  \alpha(t)\,z
+ \beta(t)\,.
\end{equation}
From \eqref{S_z+w1} for $z = 0$ it follows that
\begin{equation}\label{alpha}
\alpha(t) = - \frac{\eta_1(t)}{\omega_1}\Bigg( \sum\limits_{j_1 = 1}^{m_1} \mathcal{L}_{1,j_1}(t) + \sum\limits_{j_2 = 1}^{m_2} \mathcal{L}_{2,j_2}(t)\Bigg).
\end{equation}

Now we will demand that $\mathcal{F}(0, t)$ does not depend on $t$. It is possible since the conformal mappings of concentric annuli onto domains are defined up to rotation of the annuli. So we have  $\mathcal{F}(0, t) = 0$ and $\mathcal{H}(0, t) = 0$. Then
\begin{equation}\label{beta}
\beta(t) = \sum\limits_{j_1 = 1}^{m_1} \mathcal{L}_{1,j_1}(t)\, \zeta({z}_{1,j_1,0}(t)) + \sum\limits_{j_2 = 1}^{m_2} \mathcal{L}_{2,j_2}(t) \,\zeta({z}_{2,j_2,0}(t)).
\end{equation}
From \eqref{S_z+w2} for $z = 0$ it follows that
\begin{multline*}
\dot{\omega}_2(t) = - \alpha(t)\,\omega_2(t) - \eta_2(t)\Bigg( \sum\limits_{j_1 = 1}^{m_1} \mathcal{L}_{1,j_1}(t) + \sum\limits_{j_2 = 1}^{m_2} \mathcal{L}_{2,j_2}(t)\Bigg)\\
= \frac{\eta_1(t)\, \omega_2(t) - \eta_2(t) \, \omega_1}{\omega_1} \Bigg( \sum\limits_{j_1 = 1}^{m_1} \mathcal{L}_{1,j_1}(t) + \sum\limits_{j_2 = 1}^{m_2} \mathcal{L}_{2,j_2}(t)\Bigg),
\end{multline*}
and, with the help of Legendre's relation \eqref{zetaodnor}, we obtain
\begin{equation}\label{w2_t_it}
\dot{\omega}_2(t) = i \Bigg( \sum\limits_{j_1 = 1}^{m_1} \mathcal{L}_{1,j_1}(t) + \sum\limits_{j_2 = 1}^{m_2} \mathcal{L}_{2,j_2}(t)\Bigg).
\end{equation}
At last, from \eqref{Hab}, \eqref{alpha} and \eqref{beta} we deduce that
\begin{equation}\label{h_z_t_it}
\mathcal{H}(z, t) = \sum\limits_{j_1 = 1}^{m_1} \mathcal{L}_{1,j_1}(t)\, \mathcal{K}_{1,j_1}(z,t)\, + \sum\limits_{j_2 = 1}^{m_2} \mathcal{L}_{2,j_2}(t)\, \mathcal{K}_{2,j_2}(z,t)\,,
\end{equation}
where
$$
\mathcal{K}_{1,j_1}(z,t) = \zeta(z - {z}_{1,j_1,0}(t)) - \frac{\eta_1(t)}{\omega_1}\,z + \zeta({z}_{1,j_1,0}(t))\,,
$$
$$
\mathcal{K}_{2,j_2}(z,t) = \zeta(z - {z}_{2,j_2,0}(t)) - \frac{\eta_1(t)}{\omega_1}\,z + \zeta({z}_{2,j_2,0}(t))\,.
$$

We note that the functions  $\mathcal{K}_{1,j_1}$ and $\mathcal{K}_{2,j_2}$ are similar to the Villat kernel entering into the integral representation giving a solution to the Schwarz problem about recovering an analytic function by the boundary values of its real part (see, e.g.  \cite[ch.~10, \S.~56]{akhiezer}). Therefore, we proved the following fact.

\begin{thr}\label{family} {The smooth family of conformal mappings $\mathcal{F}(z,t)$, defined above, satisfies the PDE
\begin{equation}\label{f_p_t_f_p_z_h_it}
\frac{\dot{\mathcal{F}}(z, t)}{\mathcal{F}'(z, t)} = \mathcal{H}(z, t)\,.
\end{equation}
Here $\mathcal{H}(z, t)$ is defined by \eqref{h_z_t_it}, and the Weierstrass $\zeta$-function from \eqref{h_z_t_it} has periods  $\omega_1=2\pi$ and $\omega_2=\omega_2(t)$ satisfying \eqref{w2_t_it}.
}
\end{thr}

Now we will find a differential equation for ${z}_{1,j_1,0}(t)$. To do this, we will write $\dot{\mathcal{F}}'({z}_{1,j_1,0}(t),t)$ in two different ways. From \eqref{1l1} it follows that
$$
\mathcal{F}'(z,t) = \mathcal{F}''({z}_{1,j_1,0}(t),t)\,(z - {z}_{1,j_1,0}(t)) + \dots\,.
$$
Then
\begin{equation}\label{f_z24_t_t_pp_it}
\dot{{z}}_{1,j_1,0}(t) = - \frac{\dot{\mathcal{F}}'({z}_{1,j_1,0}(t),t)}{\mathcal{F}''({z}_{1,j_1,0}(t),t)} \,.
\end{equation}
On the other side, as $z \rightarrow {z}_{1,\ell_1,0}(t)$, by Theorem~\ref{family}, we have
\begin{equation}\label{f_z24_t_t_pp_it2}
\dot{\mathcal{F}}' ({z}_{1,\ell_1,0}(t), t) = \mathcal{F}''({z}_{1,\ell_1,0}(t), t)\,\mathcal{H}({z}_{1,\ell_1,0}(t), t) + \mathcal{F}'({z}_{1,\ell_1,0}(t), t)\,\mathcal{H}'({z}_{1,\ell_1,0}(t),
t)\,.
\end{equation}
Comparing \eqref{f_z24_t_t_pp_it} and \eqref{f_z24_t_t_pp_it2}, and taking into account \eqref{f_z_t_pp_z1k2}, we obtain
\begin{multline}\label{x_1_k_t_s}
\dot{{z}}_{1,\ell_1,0}(t) =  -\sum\limits_{j_1 = 1,\, j_1 \neq \ell_1}^{m_1} \mathcal{L}_{1,j_1}(t)\, \mathcal{K}_{1,j_1}({{z}}_{1,\ell_1,0}(t),t)\, -\sum\limits_{j_2 = 1}^{m_2}
\mathcal{L}_{2,j_2}(t)\, \mathcal{K}_{2,j_2}({{z}}_{1,\ell_1,0}(t),t)\, \\
 -\, \mathcal{L}_{1,\ell_1}(t) \Bigg[\mathcal{Q}({z}_{1,\ell_1,0}(t),t) \,+ \sum\limits_{j_1 = 1,\, j_1 \neq \ell_1}^{m_1} \mathcal{Q}_{1,i_1}({z}_{1,\ell_1,0}(t), t)\, + \sum\limits_{j_2 = 1}^{m_2}
 \mathcal{Q}_{2,i_2}({z}_{1,\ell_1,0}(t), t) \Bigg].
\end{multline}
Similarly,
\begin{multline}\label{z_2_l_t_s}
\dot{{z}}_{2,\ell_2,0}(t) =  -\sum\limits_{j_1 = 1}^{m_1} \mathcal{L}_{1,j_1}(t)\, \mathcal{K}_{1,j_1}({{z}}_{2,\ell_2,0}(t),t)\, - \sum\limits_{j_2 = 1,\, j_2 \neq \ell_2}^{m_2}
\mathcal{L}_{2,j_2}(t)\, \mathcal{K}_{2,j_2}({{z}}_{2,\ell_2,0}(t),t)\, \\
 -\, \mathcal{L}_{2,\ell_2}(t) \Bigg[\mathcal{Q}({z}_{2,\ell_2,0}(t),t) \,+ \sum\limits_{j_1 = 1}^{m_1} \mathcal{Q}_{1,i_1}({z}_{2,\ell_2,0}(t), t)\, + \sum\limits_{j_2 = 1, \, j_2 \neq \ell_2}^{m_2}
 \mathcal{Q}_{2,i_2}({z}_{2,\ell_2,0}(t), t) \Bigg].
\end{multline}

Now we will find a differential equation for $C_1(t)$. On one side, from Theorem~\ref{family} it follows that
$$
\frac{\dot{\mathcal{F}}'(0, t)}{{\mathcal{F}}'(0, t)} =  -\sum\limits_{j_1 = 1}^{m_1} \mathcal{L}_{1,j_1}(t) \, \mathcal{P}({z}_{1,j_1,0}(t), t) - \sum\limits_{j_2 = 1}^{m_2} \mathcal{L}_{2,j_2}(t)\,
\mathcal{P}({z}_{2,j_2,0}(t), t)\,,
$$
where
$\mathcal{P}(z, t) = \wp(z;\, \omega_1,\omega_2(t)) + {\eta_1(t)}/{\omega_1}$.

On the other side, from \eqref{f_sem} we have
\begin{multline*}
\frac{\dot{\mathcal{F}}'(0, t)}{{{\mathcal{F}}'(0, t)}} = \frac{\dot{C}_1(t)}{C_1(t)} \, + \sum\limits_{i_1 = 1}^{n_1} (\alpha_{1,i_1} - 1) \, \mathcal{R}(z_{1,i_1}(t),t) \, +  \sum\limits_{i_2 =
1}^{n_2} (\alpha_{2,i_2} - 1)\, \mathcal{R}(z_{2,i_2}(t),t) \, \\
+  \sum\limits_{j_1 = 1}^{m_1} \mathcal{R}({z}_{1,j_1,0}(t),t) + (\varphi_{1,j_1} - 1) \, \mathcal{R}(z_{1,j_1,1}(t),t) - \varphi_{1,j_1} \, \mathcal{R}(z_{1,j_1,2}(t),t) \,\\
+  \sum\limits_{j_2 = 1}^{m_2} \mathcal{R}({z}_{2,j_2,0}(t),t)  +  (\varphi_{2,j_2} - 1) \, \mathcal{R}(z_{2,j_2,1}(t),t)- \varphi_{2,j_2} \, \mathcal{R}(z_{2,j_2,2}(t),t)\,,
\end{multline*}
where
$$
\mathcal{R}(z(t),t) = \zeta(z(t);\,\omega_1,\omega_2(t))\, \dot{z}(t) + \dot{\omega}_2(t)\, \frac{\partial \log \sigma (z(t);\,\omega_1,\omega_2(t))}{\partial \omega_2}\,.
$$
Then
\begin{multline}\label{C_1_p}
-\frac{\dot{C}_1(t)}{C_1(t)} =  \sum\limits_{i_1 = 1}^{n_1} (\alpha_{1,i_1} - 1) \, \mathcal{R}(z_{1,i_1}(t),t)\, +  \sum\limits_{i_2 = 1}^{n_2} (\alpha_{2,i_2} - 1) \, \mathcal{R}(z_{2,i_2}(t),t) \,\\
+  \sum\limits_{j_1 = 1}^{m_1} \mathcal{R}({z}_{1,j_1,0}(t),t) + (\varphi_{1,j_1} - 1) \, \mathcal{R}(z_{1,j_1,1}(t),t)  - \varphi_{1,j_1} \, \mathcal{R}(z_{1,j_1,2}(t),t) \,\\
+  \sum\limits_{j_2 = 1}^{m_2} \mathcal{R}({z}_{2,j_2,0}(t),t)  +  (\varphi_{2,j_2} - 1) \, \mathcal{R}(z_{2,j_2,1}(t),t) - \varphi_{2,j_2} \, \mathcal{R}(z_{2,j_2,2}(t),t)\,\\
+ \sum\limits_{j_1 = 1}^{m_1} \mathcal{L}_{1,j_1}(t)\, \mathcal{P}({z}_{1,j_1,0}(t), t)\, +  \sum\limits_{j_2 = 1}^{m_2} \mathcal{L}_{2,j_2}(t)\, \mathcal{P}({z}_{2,j_2,0}(t), t)\, .
\end{multline}

Thus, we have the following statement.

\begin{thr}\label{accessory}
 {The accessory parameters of the family  $\mathcal{F}(z, t)$ satisfy the system of ODEs \eqref{z_1k_p_t_it} -- \eqref{z_2k_p_t_it1}, \eqref{w2_t_it}, \eqref{x_1_k_t_s},
 \eqref{z_2_l_t_s}  and \eqref{C_1_p} where the function $\mathcal{H}(z, t)$ is defined by \eqref{h_z_t_it}. In \eqref{C_1_p} the partial derivative of $\log \sigma(z)$ with respects to the period
 $\omega_2$ is given by the formula \eqref{spw2}.
}
\end{thr}

\begin{cor}\label{modt} {The conformal modulus $m(t) = \textup{Mod} (\mathcal{D}(t))$ satisfies the differential equation
$$
\dot{m}(t) = \frac{1}{2\omega_1} \Bigg( \sum\limits_{j_1 = 1}^{m_1} \mathcal{L}_{1,j_1}(t) + \sum\limits_{j_2 = 1}^{m_2} \mathcal{L}_{2,j_2}(t)\Bigg).
$$
}
\end{cor}

\begin{rem}\label{degen} We note that the right-hand sides of \eqref{z_1k_p_t_it} -- \eqref{z_2k_p_t_it1}, \eqref{w2_t_it}, \eqref{x_1_k_t_s}, \eqref{z_2_l_t_s}  and \eqref{C_1_p} have singularities at the initial moment (say, $t = 0$). This is caused by the fact that the lengths  of the slits tend to zero as $t\to 0$, consequently, for every fixed $1\le k\le 2$ and $1\le j_k\le m_k$ the points $z_{k,j_k,0}(t)$, $z_{k,j_k,1}(t)$ and $z_{k,j_k,2}(t)$ have the same limit. To avoid degeneration when solving the Cauchy problem for the system, in practice we change initial values of the parameters by the values of order $10^{-12}$.
\end{rem}

\begin{rem}\label{angle_basepoint}
Consider the case when we make a new $(m_1 + 1)$th slit from the $j$th vertex of $\Gamma_1$ (the case of of a slit going from a vertex of $\Gamma_2$ is studied simiklarly). Then the family of conformal mappings~\eqref{f_sem} will be rewritten in the form
\begin{multline*}
\mathcal{F}(z,t) = C_1(t) \int\limits_{0}^{z} \exp\{c(t)\, \xi \} \prod\limits_{i_1 = 1,\,i_1\neq j}^{n_1} \sigma^{\alpha_{1,i_1} -1}(\xi - z_{1,i_1}(t)) \prod\limits_{i_2 = 1}^{n_2} \sigma^{\alpha_{2,i_2} -1}
(\xi - z_{2,i_2}(t))\,\times \\
\times \prod\limits_{j_1 = 1}^{m_1 +1} s_{1,j_1}(\xi,t) \prod\limits_{j_2 = 1}^{m_2} s_{2,j_2}(\xi,t)\,  d\xi + C_2\,,
\end{multline*}
where
$$
s_{1,j_1}(z,t) = \sigma(z - {z}_{1,j_1,0}(t))\, \sigma^{\varphi_{1,j_1} -1}(z - z_{1,j_1,1}(t))\,\sigma^{-\varphi_{1,j_1}} (z - z_{1,j_1,2}(t))\, , \quad 1 \le j_1 \le m_1,
$$
$$
s_{1,m_1+1}(z,t) = \sigma(z - {z}_{1,m_1+1,0}(t))\, \sigma^{\varphi_{1,m_1+1,1} - 1}(z - z_{1,m_1+1,1}(t))\,\sigma^{\varphi_{1,m_1+1,2}-1}(z - z_{1,m_1+1,2}(t))\,  ,
$$
$$
s_{2,j_2}(z,t) = \sigma(z - {z}_{2,j_2,0}(t))\, \sigma^{\varphi_{2,j_2} - 1}(z - z_{2,j_2,1}(t))\,\sigma^{-\varphi_{2,j_2}}(z - z_{2,j_2,2}(t))\,  ,
$$
\begin{multline*}
c(t) = \frac{\eta_1(t)}{\omega_1}\Bigg[\sum\limits_{i_1=1,\, i_1\neq j}^{n_1} (\alpha_{1,i_1} -1)\,z_{1,i_1}(t)+\sum\limits_{i_2=1}^{n_2} (\alpha_{2,i_2} -1)\,x_{2,i_2}(t)\, +\\
+ \sum\limits_{j_1=1}^{m_1} \Big({z}_{1,j_1,0}(t) + (\varphi_{1,j_1}-1)\,z_{1,j_1,1}(t) -\varphi_{1,j_1} \, z_{1,j_1,2}(t)\Big) + \\
+ z_{1,m_1+1,0}(t) + (\varphi_{1,m_1+1,1}-1)z_{1,m_1+1,1}(t) + (\varphi_{1,m_1+1,2} - 1)z_{1,m_1+1,2}(t) \,+\\
+  \sum\limits_{j_2=1}^{m_2}\Big({x}_{2,j_2,0}(t) + (\varphi_{2,j_2} - 1)\, x_{2,j_2,1}(t) - \varphi_{2,j_2}\,x_{2,j_2,2}(t) \Big) \Bigg] + \eta_2(t)\,.
\end{multline*}
The interior angles of the domain which is cut along the slit  satisfy the relation $\varphi_{1,m_1+1,1} + \varphi_{1,m_1+1,2} = \alpha_{1,j}$. The ODE system for this family can be also rewritten in a suitable form. 
\end{rem}


\section{Examples}\label{sect_ex}

In this section we give some numerical examples showing the efficiency of our method. For numeric calculations we used the Wolfram Mathematica package, version 13.3.

\begin{ex}\label{ex5.1}
Consider one parameter family of conformal mappings of annuli $\mathcal{A}(t) = \{\tau: \; q(t) < |\tau| < 1\}$ onto doubly connected polygonal domains $\mathcal{D}(t)$, $0\le t\le 1$. The  exterior boundary component of every domain $\mathcal{D}(t)$ is the boundary of the fixed rectangle $[-1,1]\times [-0.5,0.5]$, and the interior one is a rectilinear slit along the segment $[a_1(t), a_2(t)]$ of the real axis (Fig.~\ref{f1}). We will assume that the segment $[a_1(0), a_2(0)]$ coincides with $[-0.5,0.5]$ and the endpoints of the slit move with constant velocities $v_1=\dot{a}_{1}(t)$ and $v_2=\dot{a}_{2}(t)$. Therefore, we have  $a_1(t)=v_1 t-0.5$, $a_2(t)=v_2 t+0.5$, $0\le t \le 1$. The family of such conformal mappings has the form
\begin{equation}\label{confP}
\mathcal{F}(z,t) = C_1(t) \int\limits_{0}^{z} \exp\{c \xi \} \prod\limits_{i_1 = 1}^{4} \sigma^{-1/2} (\xi - z_{1,i_1}(t))\, \prod\limits_{i_2 = 1}^{2} \sigma(\xi - z_{2,i_2}(t))\,  d\xi + C_2(t)\,,
\end{equation}
where $z = -i \ln \tau$ and $c(t) = \eta_1(t)/2 + \eta_2(t)$.
In~\eqref{confP} the accessory parameters $z_{1,i_1}(t)$, $1 \le i_1 \le 4$, are symmetrical with respect to the imaginary axis, i.e. $z_{1,1}(t) = - z_{1,4}(t)$ and $z_{1,2}(t) = - z_{1,3}(t)$. The parameters $z_{2,1}(t)$ and $z_{2,2}(t)$ only depend on the period $\omega_2(t)$, i.e. $z_{2,1}(t) =\omega_1/2 + \omega_2(t)/2$ and $z_{2,2}(t) = \omega_2(t)/2$.

\begin{figure}[h!]\label{f1}
\center{\includegraphics[width=0.49\linewidth]{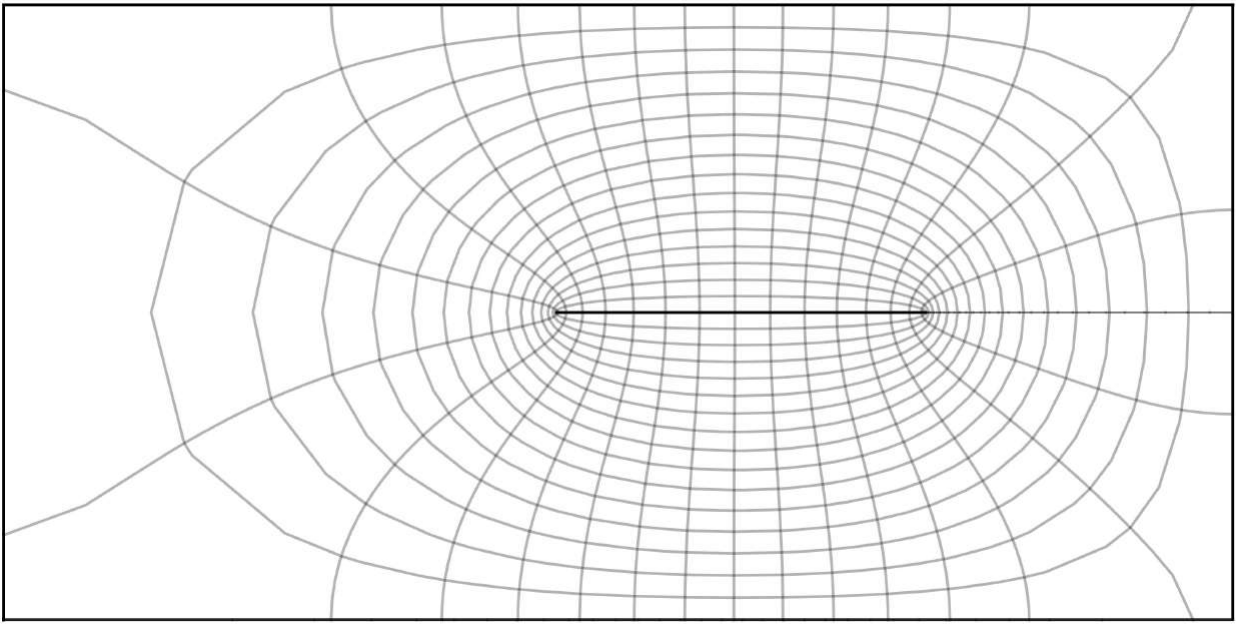}}
\begin{flushleft}
\caption{The image of the polar grid under the conformal mapping of an annulus onto the rectangle $(-1,1)\times (-0.5,0.5)$ with a slit along the segment $[-0.1, 0.5]$.}
\end{flushleft}
\end{figure}

We note that earlier in Example~\ref{rect_in}, with the help of elliptic functions,  we found the values of parameters for conformal mappings of such domains.  Here we apply the results of Section~\ref{sect_fam}  and determine the values of parameters at $t=1$, considering that the initial values of the parameters (for $t=0$) are known; we will take them from Table~\ref{tab1}.

From Theorem~\ref{accessory} and Remark~\ref{angle_basepoint} it follows that the accessory parameters satisfy the following differential equations
$$
\dot{\omega}_2(t) = i \sum\limits_{j_2 = 1}^{2} \mathcal{L}_{2,j_2}(t)\,,\quad \dot{z}_ {1,k}(t) = -\sum\limits_{j_2 = 1}^{2} \mathcal{L}_{2,j_2}(t)\, \mathcal{K}_{2,j_2} (z_{1,k}(t),t)\,, \quad 1 \le k \le 2\,,
$$
$$
\frac{\dot{C}_1(t)}{C_1(t)} =\sum\limits_{i_1 = 1}^{2} \mathcal{R}(z_{1,i_1}(t),t ) - \sum\limits_{i_2 = 1}^{2} \mathcal{R}(z_{2,i_2}(t),t) - \sum\limits_{j_2 = 1}^{2} \mathcal{ L}_{2,j_2}(t) \, \mathcal{P}({z}_{2,j_2}(t), t)\,.
$$
Solving the Cauchy problem for this ODE system, we obtain the values of the accessory parameters and conformal modulus corresponding to the value $t = 1$; they are shown in Table~\ref{t51} for different speeds of the endpoints of the slit.

\begin{table}[ht]\label{t51}
\caption{Some values of the accessory parameters and the conformal module for the conformal mapping of the annulus onto the rectangle $(-1,1) \times (-0.5,0.5)$ with a cut along the segment~$[v_1 - 0.5, v_2 + 0.5]$.}
\begin{center}
{
\scriptsize
\begin{tabular}{|{c}|{c}|{c}|{c}|{c}|{c}|{c}|}
\hline
$v_1$  & $v_2$   & $\omega_2(1)$ & $z_{1,1}(1)$ & $z_{1,2}(1)$ & $m(1)$ \\
\hline
$-0.25$& $0.25$  & $1.5034301873176772\,i$ & $0.47104824707292203$ & $2.6705444064077497$ & $0.11963917295259123$ \\
\hline
$0.25$ & $-0.25$ & $3.3301140782903405\,i$ & $0.1992807585816138$ & $2.9423118950085363$ & $0.2650020583099093$ \\
\hline
$0.4$  & $0$     & $2.9794506900288944\,i$ & $0.39772855448531813$ & $3.0299504820063854$ & $0.2370971525083285$ \\
\hline
$0.6$  & $0.4$   & $2.118321160564505\,i$ & $1.251124799561501$ & $3.099177418082569$ & $0.16857064188000076$ \\
\hline
\end{tabular}
}
\end{center}
\end{table}

Comparison of the values from Tables~\ref{tab1} and~\ref{t51} shows a good coincidence, up to 9 digits after the decimal point. Numerical calculations show that the accuracy of calculations falls when the endpoints of the slit are close enough to the boundary of the rectangle, meanwhile the conformal modulus is very close to zero.

\end{ex}

\begin{ex}\label{ex5.2}
Consider the problem of finding conformal mapping of an annulus $\mathcal{A} = \{\tau: \; q < |\tau| < 1\}$ onto the doubly connected domain $\mathcal{D}_1$ that is obtained from the interior of the rectangle $\Pi = (-1,1) \times (-0.5,0.5)$ by removing the triangle with vertices at the points
$w_{2,1} = 0.25$, $w_{2,2} = -0.25$ and $w_{2,3} = 0.25\, i$ (see Fig.~\ref{p2}).
By Theorem~\ref{int_repr}, the conformal mapping of $\mathcal{D}_1$ has the form
\begin{equation}\label{confG}
\mathcal{F}(z) = C_1 \int\limits_{0}^{z} \exp\{c \xi \} \prod\limits_{i_1 = 1}^{4} \sigma^{-1/2} (\xi - z_{1,i_1})\,  \sigma^{3/4}(\xi - z_{2,1})\,\sigma^{3/4}(\xi - z_{2,2})\, \sigma^{1/2}(\xi - z_{2,3})\, d\xi + C_2\,,
\end{equation}
where
$$
c = \frac{\eta_1}{\omega_1}\Bigg[-\frac{1}{2} \sum\limits_{i_1=1}^{4} z_{1,i_1} + \frac{3}{4}\, x_{2,1} + \frac{3}{4}\, x_{2,2}  + \frac{1}{2}\, x_{2,3} \Bigg] + \eta_2\,.
$$
Here $\mathcal{F}(z_{2,j}) = w_{2,j}$, $1 \le j \le 3$.

\begin{figure}[h!]
\center{\includegraphics[width=0.49\linewidth]{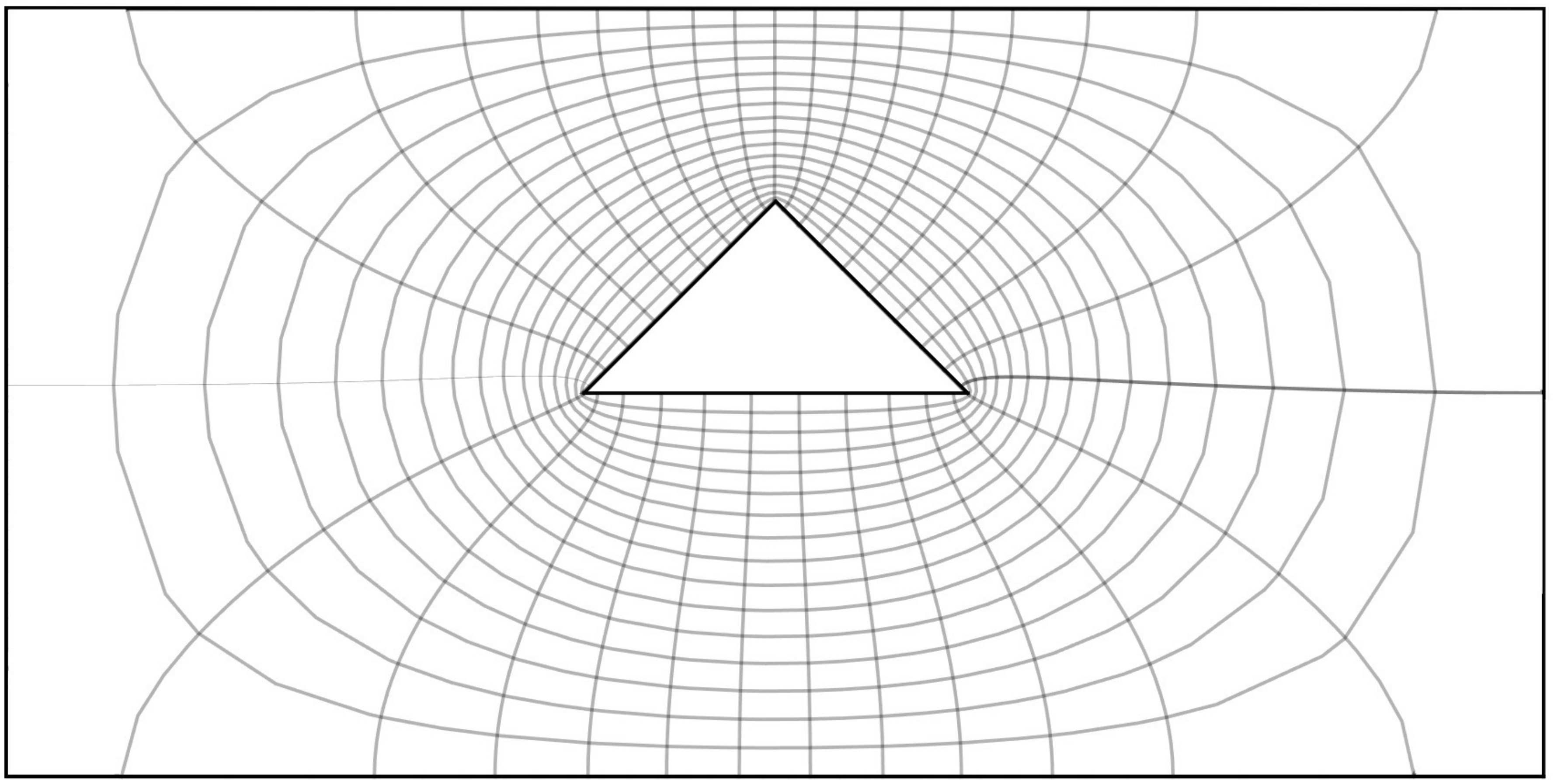}}
\caption{The image of the polar grid under the conformal mapping of the annulus onto the domain obtained from the rectangle $(-1,1) \times (-0.5,0.5)$ by removing the triangle with vertices at the points $0.25$, $-0.25$ and ~$0.25 \,i$.}
\label{p2}
\end{figure}

We will find approximate values of the parameters in \eqref{confG}.
First we consider the conformal mapping of an annulus onto the domain $\mathcal{D}_0$ that is the rectangle $\Pi$ with a rectilinear slit along the segment $[-0.25, 0.25]$; all parameter for this mapping are known (Example~\ref{rect_in}). Then we will release rectilinear slits from the points $w_{2,1}$ and $w_{2,2}$. Let the endpoints of the slits move with constant velocities according to the laws  $v_{1}(t)=0.25(1+(i-1)t)$ and  $v_2(t)=0.25(-1+(i+1)t)$. Then, at the moment $t=1$ they meet at the point $w_{2,3}$. Denote by $\mathcal{D}(t)$ the doubly connected domain that is obtained from $\Pi$ by making the slit along  the polygonal line with vertices at the points $w_{2,1}$, $w_{2,2}$, $v_1(t)$ and $v_2(t)$. We note that $\mathcal{D}(0)=\mathcal{D}_0$ and, as $t\to 1$, the family of domains $\mathcal{D}(t)$ converges to $\mathcal{D}_1$ as to a kernel in the sense of Caratheodory (see \cite[ch.~2, \S~5;  ch.~5, \S~5]{Goluzin_book}).

The family of conformal mappings of annuli $\mathcal{A}(t) = \{\tau: \; q(t) < |\tau| < 1\}$ onto $\mathcal{D}(t)$, $0<t<1$, has the form
\begin{multline}\label{confDt}
\mathcal{F}(z,t) = C_1(t) \int\limits_{0}^{z} \exp\{c(t)\, \xi \} \prod\limits_{i_1 = 1}^{4} \sigma^{-1/2} (\xi - z_{1,i_1}(t))\,  \sigma(\xi - z_{2,1,0}(t))\,\sigma^{3/4}(\xi - z_{2,1,1}(t)) \\
\times \, \sigma^{-3/4}(\xi - z_{2,1,2}(t))\,\sigma(\xi - z_{2,2,0}(t))\,\sigma^{-3/4}(\xi - z_{2,2,1}(t)) \, \sigma^{3/4}(\xi - z_{2,2,2}(t))\, d\xi + C_2(t)\,,
\end{multline}
where
\begin{multline*}
c(t) = \frac{\eta_1(t)}{\omega_1}\Bigg[-\frac{1}{2} \sum\limits_{i_1=1}^{4} z_{1,i_1}(t) + x_{2,1,0}(t) + \frac{3}{4}\, x_{2,1,1}(t) - \frac{3}{4}\, x_{2,1,2}(t) \, \\
+\, x_{2,2,0}(t) - \frac{3}{4}\, x_{2,2,1}(t) + \frac{3}{4}\, x_{2,2,2}(t) \Bigg] + \eta_2(t)\,.
\end{multline*}
Here the parameters $z_{1,i_1}(t)\; (\text{mod}\, \omega_1)$, $1 \le i_1 \le 4$, are real and taken in the ascending order on $[0,\omega_1]$, and the parameters $z_{2,k,j}(t) = x_{2,k,j}(t) + \omega_2(t)/2$, $1 \le k \le 2, 0 \le j \le 2$, are such that $x_{2,k,j}(t)\; (\text{mod}\, \omega_1)$ are real and taken in the ascending order on $[0,\omega_1]$. As we noted above, for $t=0$ the values of the parameters coincide with the corresponding values of parameters for the conformal mapping of an annulus onto $\mathcal{D}_0$; we give them from Table~\ref{tab1}, taking into account Remark~\ref{angle_basepoint}.

Solving the Cauchy problem for the system of ODEs constructed for the family~\eqref{confDt} (see Theorem~\ref{accessory} and Remark~\ref{angle_basepoint}), we obtain the needed values of the accessory parameters and the conformal modulus of $\mathcal{D}_1$ at $t=1$; they are given in Table~\ref{t52}. We also determine the values $c(1) = -0.1455922673170743 - 1.2449331197330842\, i$ and $C_1(1) = 0.1299350302810949 + 0.11027609828928062 \, i$.

\begin{table}
\caption{Values of accessory parameters and conformal module for mapping of an annulus onto a rectangle $(-1,1) \times (-0.5,0.5)$ with the triangle removed.
}
{\scriptsize
\begin{center}
\begin{tabular}{|{c}|{c}|{c}|{c}|}
\hline
$m(1)$  & $\omega_2(1)$ & $z_{1,1}(1)$ & $z_{1,2}(1)$ \\
\hline
$0.21029028897501653$ & $2.6425857078607464\, i$ & $0.185196737529923$ & $3.4688730938781616$ \\
\hline
$z_{1,3}(1)$  & $z_{1,4}(1)$ & $x_{2,2,2}(1)$ &  $x_{2,1,1}(1)$ \\
\hline
$3.8320979743892862$ & $6.105157164198392$ & $6.109642627793988$ & $3.827612510793681$ \\
\hline
$x_{2,1,0}(1)$  & $x_{2,1,2}(1)$ & $x_{2,2,1}(1)$ &  $x_{2,2,0}(1)$ \\
\hline
$1.8334660411846952$ & $1.8270349281572866$ & $1.8270349032507693$ & $1.8206037902233696$ \\
\hline
\end{tabular}
\end{center}
}\label{t52}
\end{table}

As approximate values of the accessory parameters in \eqref{confG} we can take $z_{1,i_1} = z_{1,i_1}(1)$, $1 \le i_1 \le 4$, $z_{2,1} = z_{2,2,2}(1)$, $z_{2,2} = z_{2,1,1}(1)$; for the value $z_{2,3}$ we can take either the arithmetic mean of the numbers $z_{2,1,0}(1)$, $z_{2,1,2}(1)$, $z_ {2,2,1}(1)$ and~$z_{2,2,0}(1)$ or any of these values. We note that, as~$t\to1$, the accessory parameters $z_{2,1,0}$, $z_{2,1,2}$, $z_{2,2,1}$ and~$z_{2, 2,0}$ approach each other indefinitely and merge into a unique  point, but in numerical calculations their values differ by values of the order of $10^{-4}$; this is caused by the non-uniform convergence of the family of conformal mappings near the specified point. Note that the accuracy of calculations of polygon vertices was up to ~$10^{-6}$.

In Fig.~\ref{p2} we show the polygonal domain $\mathcal{D}_1$  and the image of the polar grid under the conformal mapping \eqref{confG}.
\end{ex}

\begin{ex}\label{53}
Consider the problem of constructing a conformal mapping of an annulus onto the domain $\mathcal{D} = \Pi_1 \backslash\Pi_2$, where $\Pi_1 = (-1,1) \times (-0.5,0.5)$ and~$\Pi_2 = [0.1,0.5]\times [0.1,0.2]$ (Fig.~\ref{p6}). The solution consists of four steps. On every step we solve the Cauchy problem for a system of ODEs; the accessory parameters obtained on the previous step, except for the final one,  define the initial values of the parameters on the next step.  For convenience of presentation, we will assume that on every step the parameter $t$ is on the segment $[0,1]$.

\text{\textbf{Step~1.}}  As an initial, we take the conformal mapping of an annulus onto  the rectangle $\Pi_1$ with the slit along the segment $[-0.1, 0.5]$;  we take  the corresponding values of parameters from Table~\ref{tab1}. We consider the family of conformal mappings of annuli onto rectangle $\Pi_1$ with an elongating slit along a two-link polyline composed from the segments $[-0.1,0.5]$ and $[0.5,0.5+0.2 t\, i]$.

\text{\textbf{Step~2.}} Now we consider the family of conformal mappings of annuli onto $\Pi_1$ with the slit along the polyline consisting of three segments,
$$
[-0.1 +0.6t,0.5], \ [0.5, 0.5 + 0.2\, i]\ \text{and}\ [0.5 + 0.2\, i, 0.5-0.4 t + 0.2\, i].
$$
We note that, with the growth of the parameter $t$, the first segment shortens and shrinks to the point $0.5$ as $t\to 1 $, and the third one grows; moreover,
it shrinks to the point  $0.5+2\,i$ as $t\to 0 $.

\text{\textbf{Step~3.}}  On this step, the slit consists of three segments,
$$
[0.5+0.1t\,i, 0.5 + 0.2\, i], \ [0.5+ 0.2\, i, 0.1 + 0.2\, i] \ \text{and} \ [0.1 + 0.2\, i,0.1 +0.2\, i-0.1t\,i].
$$
With the growth of $t$, the length of the first segment  decreases from $0.2$ to $0.1$, and the length of the third one grows from $0$ to $0.1$.

\text{\textbf{Step~4.}} On the final step, the slit consists of three fixed segments $[0.5+0.1\,i, 0.5 + 0.2\, i]$, $[0.5+ 0.2\, i, 0.1 + 0.2\, i]$, $[0.1 + 0.2\, i,0.1 +0.1\, i]$ and a growing one,  $[0.1+0.1\,i, 0.1+0.1\,i+0.4 t \,i]$. As $t\to 1$, the slit tends to the boundary of $\Pi_2$ and, by the Caratheodory convergence theorem,  the limiting function maps an annulus onto the needed domain (Fig.~\ref{p6}). We also find the values
$$
c(1) = -0.000055631445885073245-1.169959854131041\, i,
$$
$$
C_1(1) = -0.04708877454319929+0.12479568610390228\, i.
$$
The family of conformal mappings has the form
\begin{multline*}\label{confP12}
\mathcal{F}(z,t) = C_1(t) \int\limits_{0}^{z} \exp\{c(t)\, \xi \} \prod\limits_{i_1 = 1}^{4} \sigma^{-1/2} (\xi - z_{1,i_1}(t))\,  \sigma(\xi - z_{2,1}(t))\,\sigma^{-1/2}(\xi - z_{2,2}(t))\\ 
\times \, \sigma^{-1/2}(\xi - z_{2,3}(t))\,\sigma^{-1/2}(\xi - z_{2,4}(t))\,\sigma(\xi - z_{2,5}(t)) \, \sigma^{1/2}(\xi - z_{2,6}(t))\,\\ 
\times \, \sigma^{1/2}(\xi - z_{2,7}(t))\,\sigma^{1/2}(\xi - z_{2,8}(t))\, d\xi + C_2(t)\,,
\end{multline*}
where
\begin{multline*}
c(t) = \frac{\eta_1(t)}{\omega_1}\Bigg[-\frac{1}{2} \sum\limits_{i_1=1}^{4} z_{1,i_1}(t) +x_{2,1}(t) - \frac{1}{2}\, x_{2,2}(t) - \frac{1}{2}\, x_{2,3}(t) \, - \frac{1}{2}\, x_{2,4}(t) \\
+\, x_{2,5}(t) \, + \frac{1}{2}\, x_{2,6}(t) \, + \frac{1}{2}\, x_{2,7}(t) \, + \frac{1}{2}\, x_{2,8}(t) \Bigg] + \eta_2(t)\,.
\end{multline*}
The values of the parameters for $\mathcal{F}(z,1)$ are given in Table~\ref{t53}.

As a result, we obtain that the conformal mapping of an annulus onto the domain~$\mathcal{D}$ has the form
$$
\mathcal{F}(z) = C_1 \int\limits_{0}^{z} \exp\{c \xi \} \prod\limits_{i_1 = 1}^{4} \sigma^{-1/2} (\xi - z_{1,i_1})   \prod\limits_{i_2 = 1}^{4} \sigma^{1/2} (\xi - z_{2,i_2}) \, d\xi + C_2\,,
$$
where
$$
c = \frac{\eta_1}{\omega_1}\Bigg[-\frac{1}{2} \sum\limits_{i_1=1}^{4} z_{1,i_1} + \frac{1}{2} \sum\limits_{i_2=2}^{4} x_{2,i_2}  \Bigg] + \eta_2\,.
$$
Here $z_{1,i_1} = z_{1,i_1}(1)$, $1 \le i_1 \le 4$, $z_{2,1} = x_{2,1}(1) + \omega_2(1)/2$, $z_{2,2} = x_{2,6}(1) + \omega_2(1)/2$, $z_{2,3} = x_{2,7}(1) + \omega_2(1)/2$ and $z_{2,4} = x_{2,8}(1) + \omega_2(1)/2$.

\begin{figure}[h!]
\center{\includegraphics[width=0.49\linewidth]{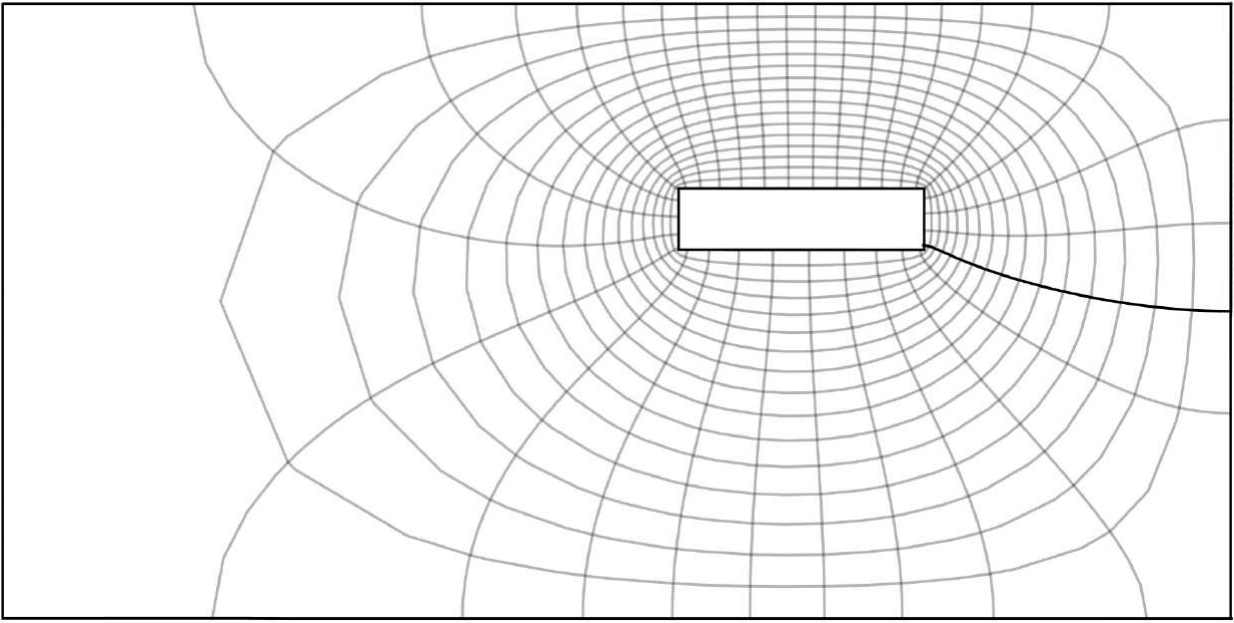}}
\caption{The image of the polar grid under the  conformal mapping of an annulus onto the  domain that is the difference of the rectangles $(-1,1) \times (-0.5,0.5)$ and~$[0.1,0.5] \times [0.1,0.2]$.}
\label{p6}
\end{figure}

\begin{table}
\caption{Values of the accessory parameters and the conformal module for mapping of an annulus onto a domain that is the difference of the rectangles $(-1,1) \times (-0.5,0.5)$ and~$[0.1,0.5] \times [0.1,0.2]$.}
{\scriptsize
\begin{center}
\begin{tabular}{|{c}|{c}|{c}|{c}|}
\hline
\multicolumn{2}{|c|}{$m(1)$} &  \multicolumn{2}{c|}{$\omega_2(1)$}  \\
\hline
\multicolumn{2}{|c|}{$0.22376354710663857$} &  \multicolumn{2}{c|}{$2.8118956629256373\, i$}  \\
\hline
$z_{1,1}(1)$  & $z_{1,2}(1)$ & $z_{1,3}(1)$ &  $z_{1,4}(1)$  \\
\hline
 $0.49730142210229983$ & $4.088563026275288$ & $4.217943920435873$ & $5.9070565839553115$ \\
\hline
$x_{2,1}(1)$  & $x_{2,2}(1)$ & $x_{2,3}(1)$ &  $x_{2,4}(1)$  \\
\hline
$6.171115702928762$ & $6.16112894071675$ &  $6.161128842028749$ &  $6.1611288420232455$ \\
 \hline
 $x_{2,5}(1)$  & $x_{2,6}(1)$ & $x_{2,7}(1)$ &  $x_{2,8}(1)$  \\
\hline
$6.157826022407041$ & $4.303856804569505$ &  $3.408905071936863$ & $0.8254470359051237$ \\
 \hline
\end{tabular}
\end{center}
}\label{t53}
\end{table}

\end{ex}

\begin{ex}\label{54}
Using the method applied in Example~\ref{53} we also calculated accessory parameters for the difference of two the rectangles, $\Pi_1 = (0,7) \times (0,4)$ and $\Pi_2 = [a,c] \times [b,d]$ for different values of $a$, $b$, $c$ and $d$. We note that such doubly connected domains $\mathcal{D}$ were considered by H.~Hakula, A.~Rasila and M.~Vuorinen~\cite{HRV}. We also calculated conformal moduli $\text{Mod} (\mathcal{D})$ and capacities $\text{Cap} (\mathcal{D})=(\text{Mod} (\mathcal{D}))^{-1}$ for such domains. In Table~\ref{t54} we give values of the moduli and capacities for some values of $a$, $b$, $c$ and $d$. Comparison of our results with those from \cite{HRV} shows a good coincidence up to $10^{-9}$. We also note that for symmetric domains the accuracy is higher, it is about $10^{-10}$.

\begin{table}\label{t54}
\caption{Some values of the conformal moduli and capacities of the domain $\mathcal{D} = \Pi_1 \backslash \Pi_2$, where~$\Pi_1 = (0.7) \times (0.4)$ and~$\Pi_2 = [a, c]\times [b,d]$.}
{\scriptsize
\begin{center}
\begin{tabular}{|{c}|{c}|{c}|{c}|{c}|{c}|{c}|{c}|{c}|}
\hline
\multicolumn{4}{|c|}{} & \multicolumn{2}{c|}{Our method} &  \multicolumn{2}{c|}{Results from~\cite{HRV}}  \\
\hline
\phantom{a}$a$\phantom{a}  & \phantom{a}$b$\phantom{a}   & \phantom{a}$c$\phantom{a} & \phantom{a}$d$\phantom{a} & $\text{Mod} (\mathcal{D})$ & $\text{Cap} (\mathcal{D})$ & $\text{Mod} (\mathcal{D})$ & $\text{Cap} (\mathcal{D})$\\
\hline
 $1$ & $1$ & $2$ & $2$ & $0.1919267753916537$ & $5.210320435798281$ & $0.19192677723893617$ & $5.210320385649294$\\
\hline
 $1$ & $1$ & $3$ & $2$ & $0.14823482589366233$ & $6.746053054478302$ & $0.1482348209832967$ & $6.746053277945276$\\
\hline
 $1$ & $1$ & $4$ & $2$ & $0.1209178457885679$ & $8.270077865500184$ & $0.12091783807693261$ & $8.27007839293125$\\
\hline
 $1$ & $1$ & $5$ & $2$ & $0.10139511296046277$ & $9.862408264093874$ & $0.10139510359023973$ & $9.86240917550835$\\
\hline
 $2$ & $1$ & $3$ & $2$ & $0.21312544007028156$ & $4.6920724230304645$ & $0.21312544403732966$ & $4.692072335693745$\\
\hline
 $2$ & $1$ & $4$ & $2$ & $0.16046010545976827$ & $6.232078666125066$ & $0.16046010434924893$ & $6.232078709256309$\\
\hline
 $2$ & $1$ & $5$ & $2$ & $0.1277611567929334$ & $7.827105085004295$ & $0.12776115200936297$ & $7.827105378062926$\\
\hline
 $3$ & $1$ & $4$ & $2$ & $0.21639757475743104$ & $4.621123878679977$ & $0.21639757713707608$ & $4.621123827863167$\\
\hline
 $3$ & $1$ & $5$ & $2$ & $0.16046010585211312$ & $6.23207865088686$ & $0.16046010434924884$ & $6.232078709256313$\\
\hline
 $4$ & $1$ & $5$ & $2$ & $0.21312544095079175$ & $4.692072403645554$ & $0.21312544403732897$ & $4.69207233569376$\\
\hline
\end{tabular}
\end{center}}
\end{table}

\end{ex}

\section{Discussion}

There are many interesting questions concerning constructing conformal mapping of polygonal domains. Our next problem is to consider one parameter families of unbounded doubly connected polygonal domains and corresponding conformal mappings. It is also interesting to  increase the accuracy of calculation of calculating accessory parameters and conformal moduli.

\section*{FUNDING}

This research is supported by the Russian Science Foundation, grant No~23-11-00066.

\end{document}